\documentclass[12pt]{amsproc}

\usepackage{amsmath}
\usepackage{amssymb}
\usepackage{amsthm}
\usepackage{amsfonts}
\usepackage{color}
\usepackage[T2A]{fontenc}
\usepackage[cp1251]{inputenc}
\usepackage[english]{babel}
\RequirePackage{srcltx}
\def\leq {\leqslant}
\def\le {\leqslant}
\def\ge {\geqslant}
\def\geq {\geqslant}

\newcommand{\sgn}{\mathop{\rm sgn}\nolimits}

\newcommand{\lef}{\Big}
\newcommand{\righ}{\Big}

\usepackage{amssymb}
\usepackage{amsfonts}
\usepackage{amsmath}
\usepackage{graphicx}%

\newtheorem{thm}{Theorem}

\newtheorem{cor}{Corollary}
\newtheorem{lem}{Lemma}
\newtheorem*{exa}{Example}
\newtheorem{rem}{Remark}

\def\cF{\mathcal{F}}

\def\cS{\mathcal{S}}

\def\N{{\mathbb N}}

\def\P{{\mathbb P}}

\def\R{{\mathbb R}}

\def\Z{{\mathbb Z}}

\def\n{\nonumber}

\def\supp{{\rm supp \, }}
\def\sign{{\rm sign \, }}

\newcommand{\bp}{{\bf p}}
\newcommand{\bq}{{\bf q}}

\begin{document}

\author{Mikhail Dyachenko}
\address{M. Dyachenko \\
Moscow State University\\ Vorobe`vy Gory, 117234, RUSSIA
and Moscow Center for Fundamental and Applied Mathematics}
\email{dyach@mail.ru}

\author{Erlan Nursultanov}
\address{E. Nursultanov \\
 Lomonosov Moscow State University (Kazakh Branch) and
Gumilyov Eurasian National University, Munatpasova 7
\\010010
Astana Kazakhstan} \email{er-nurs@yandex.ru}
\author{Sergey Tikhonov}

\address{S. Tikhonov \\
 Centre de Recerca Matem\`atica\\
Campus de Bellaterra, Edifici~C 08193 Bellaterra (Barcelona), Spain; ICREA, Pg.
Llu\'is Companys 23, 08010 Barcelona, Spain, and Universitat Aut\`onoma de
Barcelona.}
\email{stikhonov@crm.cat}

\author{Ferenc Weisz}
\address{F. Weisz \\
Department of Numerical Analysis, E\"otv\"os L. University \\
H-1117 Budapest, P\'azm\'any P. s\'et\'any 1/C., Hungary
} \email{weisz@inf.elte.hu}

  \thanks{
  This research was partially supported by
PID2020-114948GB-I00,  2017 SGR 358,
 the CERCA Programme of the Generalitat de Catalunya,
 Severo Ochoa and Mar\'{i}a de Maeztu Program for Centers and Units of Excellence in R$\&$D (CEX2020-001084-M), 
 Ministry of Education and Science of the Republic of Kazakhstan (AP051 32071, AP051 32590), and
 the Hungarian National Research, Development and Innovation Office - NKFIH, KH130426. E.N., S.T., F.W. would like to thank the Isaac Newton Institute for Mathematical Sciences, Cambridge, for support and hospitality during the programme "Approximation, sampling and compression in data science", where work on this paper was undertaken. This program was supported by EPSRC grant no EP/R014604/1.}


\title[Hardy--Littlewood-type theorems for 
  Fourier transforms in $\R^d$]
{
{
Hardy--Littlewood-type theorems for Fourier transforms in $\R^d$}
}

\date{}

\begin{abstract}
We obtain  Fourier inequalities in the weighted $L_p$ spaces for any $1<p<\infty$  involving the Hardy--Ces\`aro and Hardy--Bellman operators. We extend these results to product Hardy spaces for $p\le 1$. Moreover, boundedness of  the Hardy-Ces\`aro and Hardy-Bellman operators in various spaces (Lebesgue, Hardy, BMO) is discussed.
One of our main tools is an appropriate version of the Hardy--Littlewood--Paley inequality $
\|\widehat{f}\|_{L_{p',q}} \lesssim \left\|f\right\|_{L_{p,q}}$. 
\end{abstract}

\maketitle



\section{Introduction}



\subsection{Fourier inequalities in Lebesgue spaces}
 For the multi-dimensional Fourier transform $\widehat{f}$,
 the following Hardy-Littlewood theorem is a counterpart of the Hausdorff--Young inequality
 $	 \big\|\widehat{f}\big\|_{p'}	 \lesssim
 \left\|f\right\|_{p},$ $1 < p\le 2,$
  and reads as follows (see \cite[Th. 2.2]{sinn} and \cite[Th.2]{Nursultanov2004}):
	\begin{equation}\label{pitt-}
		\left(\int_{-\infty}^\infty \ldots\int_{-\infty}^\infty (|t_1|\ldots |t_d|)^{p-2} \left| \widehat{f}(t)\right|^{p}\,dt\right)^{1/p} \lesssim
 \left\|f\right\|_{p},\qquad 1 < p\le 2.
	\end{equation}
In particular, (\ref{pitt-}) sharpens the well-known  inequality (see, e.g., \cite[(5.19)]{ben} and \cite[p.17]{berg})
	\begin{equation}\label{pitt--}
		\left(\int_{-\infty}^\infty \ldots\int_{-\infty}^\infty |t|^{d(p-2)} \Big|\widehat{f}(t)\Big|^{p}\,dt\right)^{1/p} \lesssim\left\|f\right\|_{p},\qquad 1 < p\le 2.
	\end{equation}
Throughout the paper, by $C$ and $C_p$ we denote positive constants, which may depend on nonessential parameters and on the dimension.
As usual, $F\lesssim G$ stands for $F\le CG$. 
 If  $F\lesssim G \lesssim F$, we write $F\asymp G.$ 

As it is well known, the case $p>2$ requires special attention. 
By  $Q_N$ define a cube centered at the origin with the edge length $N$. For  $f\in L_p(\mathbb R^d)$ we define
\begin{equation}\label{F}
(\mathfrak{F}_N f)(\xi):=\int_{Q_N}f(x)e^{-i(\xi,x)}dx.
\end{equation}
As it is well known, if  $f\in L_p (\R^d)$, $1\leq p\leq 2$, the limit  $\lim\limits_{N\to +\infty}(\mathfrak{F}_N f)$
 { exists in  $L_{p'}$ and is called the Fourier transform of  $f$, where $1/p+1/p'=1$.} There are many examples in the literature showing that the Fourier transform of
  $f\in L_p (\R^d)$,  $2 < p<\infty$, is not well-defined in the usual sense; see, e.g., \cite[Ch. XVI, \S 3]{Z}.
Moreover, one can construct a Carleman-type function \cite[Ch. IV, \S 16]{bary} so that even in the case when the Fourier transform exists, neither inequality (\ref{pitt-}) nor Hausdorff--Young inequality hold; see Appendix A.

In order to derive the corresponding analogues of inequality
  (\ref{pitt-}) in the case $p\ge 2$,
 we define the Hardy-Ces\`{a}ro and Hardy-Bellman operators. First, let  $E$ consists of all
  $n$-tuples containing only  $0$ and $1.$
Set
 $
H_\varepsilon:=H_{\varepsilon_d,d}\ldots H_{\varepsilon_2,2}H_{\varepsilon_1,1},
$
where $\varepsilon \in E$ and
$$
	\left(H_{\varepsilon_i,i} f\right)(t):=\begin{cases}\frac1{t_i}\int\limits_{0}^{t_i} f(x_1,\ldots,x_i,\ldots,x_d)\,dx_i,~~\;\;\;{\text{if}} \;\;\;\;\;\;\; {\varepsilon_i} =0\\
\int\limits_{|t_i|}^{\infty} f(x_1,\ldots,\sign(t_i)x_i,\ldots,x_d)\,\frac{ dx_i}{x_i}, ~~\;\;\;{\text{if}} \;\;\;\;\;{\varepsilon_i}=1.
\end{cases}
$$

One of the main goals of the paper is to show that certain natural averages of the Fourier transforms, namely Hardy--Ces\`{a}ro and Hardy--Bellman operators, are not only well-defined for
$L_p$-functions with any $1<p<\infty$ but the  inequalities corresponding to
estimates (\ref{pitt-}) and (\ref{pitt--}) hold true.
We formulate our main result as a Pitt-type inequality with power weights.
Let us recall the known Pitt inequality  
 (\cite[Th. 2.2]{sinn}; see also \cite{ben,DeCarli2017}). 

For  $1< r\leq q<\infty$ { and $f \in \bigcup_{1 \leq s \leq 2}L_s$,}
\begin{multline}\label{vsp11-old}
\left(\int_{-\infty}^\infty \ldots\int_{-\infty}^\infty \left((|t_1| \ldots|t_d|)^\beta \left|\widehat{f}(t)\right|\right)^q\,dt\right)^{1/q} \\
\lesssim\left(\int_{-\infty}^\infty \ldots\int_{-\infty}^\infty \Big((|x_1| \ldots|x_d|)^\alpha \left|f(x)\right|\Big)^r\,dx\right)^{1/r}
\end{multline}
provided that \begin{equation}\label{con-a1}
0 \leq  \alpha=\frac1{r'}-\frac1q-\beta<\frac1{r'}
\end{equation}
 and, additionally,  \begin{equation}\label{con-a13}
\beta \leq  0,
   \end{equation}
 which is equivalent to $   \frac1{r'}-\frac1q\le\alpha$.
For optimality of these conditions see \cite[Th. 2.2]{sinn}, \cite[Sect. 4]{benedetto}, \cite[Sects. 2,3]{dyachenko}.

{ Now we point out three well-known special cases. For $\alpha=\beta=0$, $q=r'$, $1<r \leq 2$,
\begin{align}\label{e33}
\big\|\widehat{f}\big\|_{r'}\lesssim \left\|{f}\right\|_{r}
\end{align}
which is the Hausdorff-Young inequality; for $\alpha=0$, $\beta=1-2/r$, $1<r=q \leq 2$, we obtain \eqref{pitt-}, and for $\beta=0$, $\alpha=1-2/r$, $2 \leq r=q < \infty$, we establish the dual to \eqref{pitt-}, that is,
\begin{align}\label{e17}
\big\|\widehat{f}\big\|_r\lesssim\left(\int_{-\infty}^\infty \ldots\int_{-\infty}^\infty (|x_1| \ldots|x_d|)^{r-2} \left|f(x)\right|^r\,dx\right)^{1/r}.
\end{align}
The latter inequality was also proved in the recent paper \cite{ryd} provided that $f\in L_1$.

In the next theorem, we extend inequality \eqref{vsp11-old} not assuming  the condition $\beta \leq 0$.

 }

\begin{thm}\label{t2}
	Let $1< r\leq q<\infty$ and condition (\ref{con-a1}) hold.
 If a measurable function $f $ is such that
\begin{equation}\label{a1}
\left(\int_{-\infty}^\infty \ldots\int_{-\infty}^\infty \Big((|x_1| \ldots|x_d|)^\alpha |f(x)|\Big)^r\,dx\right)^{1/r}<\infty,
\end{equation}
then, for any  $\varepsilon\in E, \; t\in\mathbb R^d$, the limit
\begin{equation*}
T_\varepsilon {f}(t):=\lim_{N\to+\infty}(H_\varepsilon \mathfrak{F}_N{f})(t)
\end{equation*}
does exist.
Moreover,
for any $\varepsilon\in E$,
\begin{multline}\label{vsp11-}
\left(\int_{-\infty}^\infty \ldots\int_{-\infty}^\infty \Big((|t_1| \ldots|t_d|)^\beta \left|T_\varepsilon{f}(t)\right|\Big)^q\,dt\right)^{1/q} \\
\lesssim\left(\int_{-\infty}^\infty \ldots\int_{-\infty}^\infty \lef((|x_1| \ldots|x_d|)^\alpha \left|f(x)\right|\righ)^r\,dx\right)^{1/r}.
\end{multline}
In particular, { if $\alpha=0$, $\beta=1-2/r$, $1<r=q<\infty$, then }
\begin{equation}\label{vsp112}
	\left(\int_{-\infty}^\infty\ldots \int_{-\infty}^\infty (|t_1|\ldots |t_d|)^{r-2} \left|T_\varepsilon f(t)\right|^{r}\,dt\right)^{1/r}	\lesssim\|{f} \|_r  .
	\end{equation}

\end{thm}

A counterpart of Theorem \ref{t2} -- now  the condition $0\leq \alpha$ is not assumed but the condition $\beta\leq 0$ is fulfilled -- is given in the next theorem.

{

\begin{thm}\label{t3}
Let $1< r\leq q<\infty$ and
\begin{align}\label{e98}
	\frac1{r'}-\frac1q\le\alpha=\frac1{r'}-\frac1q-\beta<\frac1{r'}.
\end{align}
Suppose that
\begin{equation}\label{vsp}
\int_{-\infty}^\infty \ldots\int_{-\infty}^\infty \lef((|x_1| \ldots|x_d|)^\alpha \left|f(x)\right|\righ)^r\,dx<\infty.
\end{equation}
Then for any $\varepsilon\in E$ the sequence  $(H_\varepsilon \mathfrak{F}_Nf)_N$ converges to $T_{\varepsilon}f$ in the weighted Lebesgue norm, i.e.,
	\begin{align*}
		\lim_{n\to \infty} \left(\int_{-\infty}^\infty \ldots\int_{-\infty}^\infty \lef((|t_1| \ldots|t_d|)^\beta \left|T_{\varepsilon}f(t)-H_\varepsilon \mathfrak{F}_Nf(t)\right|\righ)^q\,dt\right)^{1/q} =0.
	\end{align*}
	Moreover,	
\begin{multline}\label{vsp11---+}
\left(\int_{-\infty}^\infty \ldots\int_{-\infty}^\infty \lef((|t_1| \ldots|t_d|)^\beta \left|T_{\varepsilon}f(t)\right|\righ)^q\,dt\right)^{1/q} \\
\lesssim\left(\int_{-\infty}^\infty \ldots\int_{-\infty}^\infty \lef((|x_1| \ldots|x_d|)^\alpha \left|f(x)\right|\righ)^r\,dx\right)^{1/r}.
\end{multline}
In particular, if $\beta=0$, $\alpha=1-2/r$ and $1<r=q < \infty$, then
the sequence $(H_\varepsilon \mathfrak{F}_Nf)_N$ converges to $T_{\varepsilon}f$ in the $L_r$ norm and
\begin{equation}\label{vsp113}
		\left\|T_\varepsilon{f} \right\|_r \lesssim \left(\int_{-\infty}^\infty\ldots \int_{-\infty}^\infty (|x_1|\ldots |x_d|)^{r-2} \left|f(x)\right|^{r}\,dx\right)^{1/r}.
	\end{equation}
\end{thm}

}

\vskip 0.3cm
\begin{rem}\label{r1}
 (i) {If $\alpha<1/r'$ as in \eqref{con-a1} and \eqref{e98}, then H\"older's inequality and \eqref{a1} (respectively, \eqref{vsp}) imply that $f \in L_s^{loc}$ for all $1 \leq s<r$. Thus $H_\varepsilon \mathfrak{F}_Nf$ is well defined. }
\\
(ii)
{ Note that the boundedness of Hardy's operator in $L_p$ (see \eqref{e11} below) and Hausdorff-Young inequality \eqref{e33} imply that
\begin{align}\label{e47}
	T_\varepsilon{f}=H_\varepsilon\widehat{f} \qquad \mbox{if  $f\in L_p$, $1<p \leq 2$, or $f\in \mathcal{S}$},
\end{align}
where $\mathcal{S}$ denotes the Schwartz class $\mathcal{S}(\R^d)$. So $T_{\varepsilon}$ is an extension of the operator $f \mapsto H_\varepsilon\widehat{f}$.
Under the conditions  \eqref{con-a1}, $\beta \leq 0$, and \eqref{e47}, we always see that inequality (\ref{vsp11-}) is weaker than  (\ref{vsp11-old}). } This follows from
 Hardy's inequality for averages (see, e.g., \cite{persson})
$$
\int_{-\infty}^\infty |t|^{\beta q} \left(\frac1{|t|}
\int_{0}^t g(s)ds\right)^{q}\,dt
+
\int_{-\infty}^\infty |t|^{\beta q} \left(
\int_{|t|}^\infty
\frac{g(s)}{s}
ds\right)^{q}\,dt
 \lesssim \int_{-\infty}^\infty |t|^{\beta q}
 g(t)^{q}\,dt,
$$
where $g\ge 0$ and
   $-\frac1q<\beta<1-\frac1q$, $q>1$.
\\
(iii)
 {On the other hand,
for $\beta>0$,
inequality \eqref{vsp11-} is an extension of \eqref{vsp11-old}. Moreover, (\ref{vsp112}) extends (\ref{pitt-}) for the case $p >2$ while inequality (\ref{vsp113}) extends \eqref{e17} for the case $1<p<2$ and for wider function spaces. } Inequality (\ref{vsp113})  has been obtained  in the one-dimensional case 
 in \cite{Dyachenko2018}.
\\
(iv) We note that if $1<r \leq 2$ and \eqref{vsp} holds, then $f \in L_r^{loc}$.
Putting $\beta=0$, $\alpha=1-2/r$ and $q=r$ into \eqref{vsp11-}, we obtain \eqref{vsp113} for $2 \leq r<\infty$. On the other hand, putting $\alpha=0$, $\beta=1-2/r$, $r=q$ into \eqref{vsp11---+}, we obtain \eqref{vsp112} for $1<r \leq 2$.

\end{rem}
\vskip 0.3cm
\subsection{Fourier inequalities in Hardy spaces}


The generalization of \eqref{pitt--} to Hardy spaces is  known (see Taibleson and Weiss \cite{Taibleson1980} and Garcia-Cuerva and Rubio de Francia \cite[Corollary 7.23]{Garcia-Cuerva1985} and also Bownik and Wang \cite{Bownik2013}):
\begin{equation}\label{e1}
		\left(\int_{-\infty}^\infty\ldots \int_{-\infty}^\infty |t|^{d(p-2)} \Big|\widehat{f}(t)\Big|^{p}\,dt\right)^{1/p} \lesssim  \left\|f\right\|_{H_p(\R^{d})},\qquad 0 < p\le 1.
	\end{equation}
	The case $p=1$ is called Hardy's inequality.

Our  second goal  is to generalize \eqref{pitt-} to Hardy spaces. Here we consider
the so called product Hardy spaces $H_p=H_p(\R\times\cdots\times\R)$, that is different from $H_p(\R^{n})$. See, e.g., \cite{wk4}.




\begin{thm}\label{t7}
	If $0<p \leq 1$ and $f \in H_p$, then
\begin{equation}\label{hardy1}
		\left(\int_{-\infty}^\infty \ldots \int_{-\infty}^\infty (|t_1| \cdots |t_d|)^{p-2} \left|\widehat{f}(t)\right|^{p}\,dt\right)^{1/p} \lesssim \left\|f\right\|_{H_p}.
	\end{equation}
\end{thm}
Since $\widehat{f}$ is a locally integrable function if $f \in H_p$ with $0<p \leq 1$, the integral in Theorem \ref{t7} is well defined. Note that $H_p$ is equivalent to $L_p$ if $1<p<\infty$.
We also extend the inequality of Theorem \ref{t2} to Hardy spaces.

\begin{thm}\label{t6}
	If $0<p \leq 1$ and $f \in H_p\bigcap \bigcup_{1 \leq q \leq 2}L_q$, then
	\begin{equation}\label{hardy2}
		\left(\int_{-\infty}^\infty \ldots \int_{-\infty}^\infty (|t_1| \cdots |t_d|)^{p-2} \left|H \widehat{f}(t)\right|^{p}\,dt\right)^{1/p} \lesssim \left\|f\right\|_{H_p}.
	\end{equation}
\end{thm}

\begin{rem} (i) For $p=1$, Theorem \ref{t6} holds for all $f \in H_1$. \\
(ii) Comparing the left-hand sides in  \eqref{hardy1} and   \eqref{hardy2}, we first recall
the classical reverse Hardy inequality (see \cite[Theorem 347]{hardy}):
 for any non-negative  $g$,
\begin{multline}\label{hardy3}
	\int_{-\infty}^\infty \ldots \int_{-\infty}^\infty (|t_1| \cdots |t_d|)^{p-2} g^{p}(t)\,dt\\ \lesssim \int_{-\infty}^\infty \ldots \int_{-\infty}^\infty (|t_1| \cdots |t_d|)^{p-2} \left(H g(t)\right)^{p}\,dt,\qquad 0<p \leq 1.
\end{multline}
 Moreover, the direct Hardy inequality, which is  the reverse  inequality to \eqref{hardy3}, holds for $p=1$ (see \cite[Theorem 330]{hardy}) but  does not hold in general for $0<p<1$. 
 Consider
 $$g(x)=\left\{
                           \begin{array}{ll}
                             a_n, & \hbox{$b_n<x<b_n+d_n$,} \\
                             0, & \hbox{otherwise.}
                           \end{array}
                         \right.
$$
  Then
$$I_1= \int_{0}^\infty t^{p-2} \left|g(t)\right|^{p}\,dt=a_n^p \int_{b_n}^{b_n+d_n} t^{p-2} \,dt$$
and
$$I_2= \int_{0}^\infty t^{p-2} \left|Hg(t)\right|^{p}\,dt\ge
\int_{b_n+d_n}^\infty t^{-2} \,dt \left(\int_0^{b_n+d_n} g(z)dz\right)^{p}\,dt
= \frac{ a_n^p d_n^p}{b_n+d_n}.$$
Letting now
$b_n+d_n\asymp b_n \nearrow\infty$,\; $\frac{ b_n}{d_n}\nearrow\infty,$
and
$$
a_n^p=\left(\int_{b_n}^{b_n+d_n} t^{p-2} \,dt\right)^{-1},
$$
we arrive at
$a_n^p\asymp \frac{b_n^{2-p}}{d_n}$
and
$I_1\asymp
1$ but $I_2 
 \asymp  (\frac{ b_n}{d_n})^{1-p}\nearrow\infty.
$

Thus, in the case $\widehat{f}$ is nonnegative,  \eqref{hardy2} yields  \eqref{hardy1}  for $0<p<1$ { while for $p=1$ they are equivalent.}
Without this condition on the Fourier transform, the left-hand sides of \eqref{hardy1} and  \eqref{hardy2} are not comparable. \\
(iii) It is easy to see that inequality \eqref{hardy3} is no longer true without the condition $g\ge 0.$ For example, let
$$g(x)=\left\{
                           \begin{array}{ll}
0, & \hbox{$0<x<1$;} \\
                             a_n, & \hbox{$2n-1<x<2n$;} \\
                             -a_n, & \hbox{$2n<x<2n+1$.}
                           \end{array}
                         \right.
$$
Then
$$I_1= \int_{0}^\infty t^{p-2} \left|g(t)\right|^{p}\,dt\asymp\sum_n a_n^p n^{p-2} $$
but
$$I_2= \int_{0}^\infty t^{p-2} \left|Hg(t)\right|^{p}\,dt
\asymp\sum_n a_n^p n^{-2}.$$
Taking now $a_n=n^{(1-p)/p}$ for $n<N$ and $a_n=0$ otherwise, we have
$I_2\asymp 1$ and $I_1\asymp \ln N$.

\end{rem}

 Theorem \ref{t7} is known in the one-dimensional case, cf.  \eqref{e1}
 while Theorem \ref{t6} is new even in the one-dimensional case.
 { We note that the proof of Theorem \ref{t7} was sketched in \cite{Jawerth1986}. However, since we believe it contains some 
  gaps, we present it in Section \ref{s6}.}
The dual results to Theorems \ref{t7} and \ref{t6} are also proved, see Corollary \ref{c1} in Section \ref{s6}.

\subsection{Hardy's inequalities for averages }


Let us recall the classical Hardy's inequalities. First,
	\begin{equation}\label{e11}
		\|H_\varepsilon f\|_p \lesssim\left\|f\right\|_p, 
\quad 1<p<\infty.
	\end{equation}
If all $\varepsilon_i=0$ or all $\varepsilon_i=1$, then even more is true. Let
$$
	Hf:=H_{0} f=\frac1{t_1 \cdots t_d} \int\limits_{0}^{t_1} \ldots \int\limits_{0}^{t_d} f(x_1,\ldots,x_d)\,dx_1\ldots dx_d
$$
be the Hardy--Ces\`{a}ro operator and
$$
	Bf:=H_{1} f= \int\limits_{|t_1|}^{\infty} \ldots \int\limits_{|t_d|}^{\infty} f(\sign(t_1) x_1,\ldots, \sign(t_d) x_d)\,\frac{ dx_1}{x_1} \ldots \frac{ dx_d}{x_d}
$$
be the Hardy--Bellman operator. Then
	\begin{equation}\label{e14}
		\|Hf\|_p \lesssim\left\|f\right\|_p\quad (1<p\le\infty)\quad
\mbox{and}\quad \left\|Bf\right\|_p \lesssim \left\|f\right\|_p\quad (1\le p<\infty).
	\end{equation}
The inequalities in \eqref{e11} and \eqref{e14} can be proved by iteration using the corresponding one-dimensional results.

\begin{rem}
{ If $\widehat{f} \in L_r$ for some $1<r \leq 2$ (say $f \in \mathcal{S}$),  } then the following estimates sharpen  (\ref{e11}) for $1<p \leq 2$: 
$$\|H_\varepsilon f\|_p \lesssim\left(\int_{-\infty}^\infty\ldots \int_{-\infty}^\infty (|t_1|\ldots |t_d|)^{p-2} \left|\widehat{f}(t)\right|^{p}\,dt\right)^{1/p}\lesssim \|f\|_p.
$$
Here the first inequality follows from (\ref{vsp113}) and the second one from (\ref{pitt-}).
In particular, if
\begin{align}\label{e91}
	\big|f(x)\big|\lesssim \sum_{\varepsilon\in E} \big| H_\varepsilon f(x)\big|,\qquad x_i\ne 0
\end{align}
or, more generally, if for some $1<p \leq 2$
\begin{align}\label{e92}
	\left((|x_1|\ldots |x_d|)^{-1} \int_{x}^{2x}\ldots \int_{x}^{2x} \left|{f}(t)\right|^p\,dt\right)^{1/p}
\lesssim \sum_{\varepsilon\in E} \big| H_\varepsilon f(x)\big|,\qquad x_i\ne 0,
\end{align}
{ then
$$
\|f\|_p
\asymp \left(\int_{-\infty}^\infty\ldots \int_{-\infty}^\infty (|t_1|\ldots |t_d|)^{p-2} \left|\widehat{f}(t)\right|^{p}\,dt\right)^{1/p}.
$$
Note that that the integral of the $p$-power of the left hand side of \eqref{e92} is $\|f\|_p^p$. For example, when $f \geq 0$ is non-increasing and even in each direction, then both \eqref{e91} and \eqref{e92} hold. } Equivalences of this type are usually called Hardy--Littlewood theorems or Boas-type results and
they have been previously obtained  for monotone or general monotone functions; see \cite{jam}.

\end{rem}

Our third aim in this paper is to investigate Hardy's inequalities \eqref{e14}  for limiting cases.
First, it is clear that
\[
	\left\|Hf\right\|_{BMO} \leq C \left\|Hf\right\|_\infty \leq C \left\|f\right\|_\infty.
\]
We also  note that the operator $H$ is bounded in BMO, i.e.,
\begin{equation}\label{e22}
	\left\|Hf\right\|_{BMO} \leq C \left\|f\right\|_{BMO}, \qquad f \in BMO\bigcap \bigcup_{1<q \leq \infty}L_q.
\end{equation}
See \cite{Golubov1997,Korenovskii2007} for one-dimensional functions. In the multivariate case, this  follows from \eqref{e19} below by duality, cf. Corollary \ref{c7}.

However,  the expected fact that the operator $H$ is bounded in $H_p$ for $p\le 1$ is not true. For $p=1$ it is enough to consider the function
\[
	a(x)= \left\{
           \begin{array}{ll}
             \frac{1}{2}, & \hbox{if $0 \leq  x \leq 1$,} \\
             -\frac{1}{2},\quad & \hbox{if $1< x \leq 2$,} \\
             0, & \hbox{if $x \notin [0,2]$,}
           \end{array}
         \right.
\]
which is  an atom (see Section \ref{s5}) and note that
$\int_{\R} Ha(x)\, dx = \ln 2 \neq 0$. Hence $Ha \notin H_1$ (see \cite{Golubov1997}).
From this, $H$ is not bounded in $H_p$,  $p<1$, since otherwise, by \eqref{e14} and interpolation, we would obtain boundedness in $H_1$.

Even though $H$ is not bounded in $H_p$,
we derive the following weaker result,  which can be considered as a  generalization of inequality \eqref{e14}.

\begin{thm}\label{t8}
	If $0<p \leq 1$ and $f \in H_p\bigcap \bigcup_{1<q \leq 2}L_q$, then
	$$
		\left\|Hf\right\|_p \lesssim\left\|f\right\|_{H_p}.
	$$
\end{thm}

For the operator $B$ the situation is symmetric. It is easy to see that
 $B$ is trivially bounded from $H_1$ to $L_1$ since
$$	\left\|Bf\right\|_1 \leq C \left\|f\right\|_1 \leq C \left\|f\right\|_{H_1}.
$$
In fact, it turns out that the operator $B$ is  bounded in $H_p$, i.e.,
\begin{equation}\label{e19}
	\left\|Bf\right\|_{H_p} \lesssim \left\|f\right\|_{H_p}, \qquad 0<p \leq 1, f \in H_p\bigcap \bigcup_{1 \leq q \leq 2}L_q
\end{equation}
(see \cite{Giang1995,Golubov1997,Liflyand2009,Liflyand2019} for one-dimensional functions and \cite{Giang1997,whaus} for the multivariate case) 
 but it is not bounded from $BMO$ to $BMO$ \cite{Golubov1997}. We obtain the following weaker estimate
	$$
		\left\|Bf\right\|_{BMO} \leq C \left\|f\right\|_{\infty}, \qquad f \in L_\infty\bigcap \bigcup_{1 \leq q<\infty}L_q;
	$$
see Corollary \ref{c3} below.


\subsection{Structure of the paper }

The paper is organi\-zed as follows. In Section 2 we introduce the Lorentz and net spaces.
Section 3 contains the extension of the celebrated  Hardy--Littlewood--Paley inequality 
\begin{equation*}
\big\|\widehat{f}\,\big\|_{L_{p',q}} \lesssim \left\|f\right\|_{L_{p,q}},\quad 1<p<2,\quad0<q \le \infty,
\end{equation*}
to the range $1<p<\infty$ with the help of the net spaces (see Theorem \ref{t70}).
Section \ref{s4} is devoted to the proofs of Theorems \ref{t2} and \ref{t3}.
In Sections \ref{s5} and \ref{s6}
we discuss the needed properties of the product Hardy spaces and prove Theorems  \ref{t7}--\ref{t8}, correspondingly.

\bigskip
\section{Lorentz  and net spaces}
\label{s2}

The $L_p$ space is equipped with the quasi-norm
$$
\|f\|_p:=\left(\int_{\R^d}|f(x)|^p \, dx\right)^{1/p}, \qquad 0<p<\infty,
$$
with the usual modification for $p=\infty$. Here we integrate with respect to the Lebesgue measure.
For $n$ quasi-normed spaces $X_1,\ldots ,X_d$ of one-dimensional functions, let us denote by $(X_1,\ldots ,X_d)$ the space consisting of $n$-dimensional measurable functions for which
\[
	\left\|f\right\|_{(X_1,\ldots,X_d)} := \left\|\ldots \left\|f\right\|_{X_1}\ldots\right\|_{X_d} < \infty,
\]
where the $X_j$ 
  norm is taken with respect the $j$-th variable.

The non-increasing rearrangement of a   one-dimensional measurable function $f$ is given  by
$$
f^{*}(t) := \inf \left\{\rho: |\{|f|>\rho\}| \leq t \right\}.
$$
{ For a multi-dimensional measurable function and for fixed variables $y_1,\ldots,y_{i-1},y_{i+1},\ldots,y_d$, by $f^{*_i}(y_1,\ldots,y_{i-1},\cdot,y_{i+1},\ldots,y_d)$, we denote the non-increasing rearrangement with respect to the $j$-th variable $(j=1,\ldots,n)$. Applying the non-increasing rearrangement in all variables  consecutively, we obtain
\[
	f^{*_1,\ldots,*_d}:= \left(\left(f^{*_1}\right)^{*_2}\ldots\right)^{*_d}.
\]
Various function spaces defined with the help of
iterative rearrangements
 were considered  in many papers, see e.g.  \cite{kol, Barza2006, Barza2004, Bl, Nursultanov2000}. 

 }

 Let ${\bf p}=(p_1,\ldots,p_d)$ and ${\bf q}=(q_1,\ldots,q_d)$ with $0<p_j<\infty$ and $0<q_j \leq \infty$, $j=1,2\ldots n$. In this case we will write  $0<\bp<\infty$ and $0<\bq \leq \infty$.
The Lorentz space $L_{\bp,{\bf q}}(\R^d)$ consists of all measurable functions $f$ for which
$$
\| f \|_{L_{\bf p,q}} := \left(\int_0^\infty\ldots \left(\int_0^\infty
\left(t_1^\frac1{p_1}\ldots
t_2^\frac1{p_d}{f^{*_1\ldots*_d}}(t_1,\ldots,t_d)\right)^{q_1}\frac{dt_1}{t_1}\right)^\frac{q_2}{q_1}\ldots
\frac{dt_d}{t_d}\right)^\frac1{q_d}<\infty,
$$
where in the case  $q_i=\infty$ the integral $\left(\int_0^\infty|g(t_i)|^{q_i}\frac{dt_i}{t_i}\right)^\frac1{q_i}$ is understood as $\sup_{t_i>0}|g(t_i)|$.

This definition is  different (see \cite{yatsenko, Nursultanov2000}) from the usual definition of Lorentz spaces.
Note also that the space $L_{\bf p,q}$ with $\bp = \bq$ does not coincide
with the mixed Lebesgue space $(L_{p_1},...,L_{p_d})$. However, if $p_i=q_i=p , \; i=1,2,...,n$, then $L_{\bf p,q}=L_{p}$.

We now define the net spaces \cite{jga} (see also \cite{Nursultanov1998, Nursultanov2000}). Let us denote by $M$ the collection of all rectangles $I=I_1 \times...\times I_d$ of positive measure with sides parallel to the axes. For a 
 measurable function $f$ defined on $\R^d$, we define the average function 
  by
\[
\overline{f}(t_1,...,t_d) :=	\overline{f}(t_1,...,t_d;M) := \sup_{I \in M, |I_i| \geq t_i} \frac{1}{|I_1|\ldots |I_d|} \left| \int_{I} f(x) \, dx \right|.
\]
A measurable function belongs to the net space $N_{\bp,\bq}(M)$ if
$$
\| f \|_{N_{\bf p,q}} = \left(\int_0^\infty\ldots \left(\int_0^\infty
\left(t_1^\frac1{p_1}\ldots
t_d^\frac1{p_d}\overline{f}(t_1,\ldots,t_d;M)\right)^{q_1}\frac{dt_1}{t_1}\right)^\frac{q_2}{q_1}\ldots
\frac{dt_d}{t_d}\right)^\frac1{q_d}<\infty
$$
for $0<\bp,\bq\leq\infty$. Moreover,
$N_{\bf p,q}$ is a normed linear space.

For $p_i=p,\; q_i=q, \; i=1,2,...,n,$ we also use  the notation
$$
\| f \|_{N_{p,q}} = \left(\int_0^\infty\ldots \int_0^\infty \left( (t_1\ldots t_d)^\frac1{p} \overline{f}(t_1,\ldots,t_d;M) \right)^{q}\frac{dt_1}{t_1} \ldots\frac{dt_d}{t_d}\right)^\frac1{q}
$$
and, similarly,
$$
\| f \|_{L_{p,q}} = \left(\int_0^\infty \ldots\int_0^\infty \left( (t_1\ldots t_d)^\frac1{p}{f^{*_1...*_d}}(t_1,\ldots,t_d)\right)^{q}\frac{dt_1}{t_1}\ldots \frac{dt_d}{t_d}\right)^\frac1{q}.
$$
The next result  follows easily  from the monotonicity of $\overline{f}(t_1,...,t_d;M)$ and the following
Hardy's inequality: 
$\sum_{k=1}^\infty 2^{k\alpha}\left(\sum_{m=k}^\infty a_m \right)^q\asymp
\sum_{k=1}^\infty 2^{k\alpha}a_k^q$ with $a_k\ge 0$ and $\alpha,q>0$.
\begin{lem} \label{l0} Let  $0<p<\infty, \; 0<q\leq \infty$, then 
\begin{eqnarray*}
\| f \|_{N_{p,q}}&\asymp&\left(\sum_{k\in\mathbb Z^d}\left(2^{\frac{k_1+...k_d}p}\overline{f}(2^{k_1},...,2^{k_d};M)\right)^q\right)^{\frac 1q}
\\
&\asymp&\left(\sum_{k\in\mathbb Z^d}\left(2^{\frac{k_1+...k_d}p}\sum_{m_d=k_d}^\infty\ldots\sum_{m_1=k_1}^\infty \overline{f}(2^{m_1},...,2^{m_d};M)\right)^q\right)^{\frac 1q}.
\end{eqnarray*}

\end{lem}


\bigskip
\section{Hardy--Littlewood--Paley inequality for $1<p<\infty$} 
\label{sec:hardy_and_bellman_operators}

We will need the following interpolation theorem, which is based on Theorems 1 and 2 in \cite{Nursultanov2000}.

\begin{lem}
\label{t1-} Let $0<{\bf p}_0=(p_1^0,...,p_d^0)<{\bf p_1}=(p_1^1,...,p_d^1)<\infty$, $0<{\bf q}_0=(q_1^0,...,q_d^0),{\bf q_1}=(q_1^1,...,q_d^1)<\infty, q_i^0\neq q_i^1, \;\;i=1,...,n$, and $0<{\bf r}=(r_1,...,r_d)\leq \infty$. Suppose $\bp_\varepsilon =(p^{\varepsilon_1}_1,...,p^{\varepsilon_d}_d),\; \bq_\varepsilon =(q^{\varepsilon_1}_1,...,q^{\varepsilon_d}_d)$. If $T$ is a linear operator such that for all $\varepsilon\in E$,
$$
T:(L_{p^{\varepsilon_1}_1,1},...,L_{p^{\varepsilon_d}_d,1})\;\rightarrow\; N_{\bq_\varepsilon,\bf\infty},
$$
then
\begin{equation}\label{e13-}
T:L_{\bf p,r}\to { N}_{\bf q,r},	\end{equation}
where ${\bar \theta}=(\theta_1,...,\theta_d)\in (0,1)^d$ and
$$
\frac1{\bf p}=\frac{1-\bar\theta}{\bf {\bf p_0}}+\frac{\bar\theta}{\bf p_1},
\quad \frac1{\bf q}=\frac{1-{\bar\theta}}{\bf {\bf {\bf
q_0}}}+\frac{\bar\theta}{\bf q_1}.
$$
\end{lem}
\begin{proof}
By Lemma 4 (c) in \cite{Nursultanov2000}, we have that
$$
T:\left({\bf A_0,A_1}\right)_{{\bar \theta},\bf r}\to \left({\bf B_0,B_1}\right)_{\bf {\bar \theta},r},
$$
where
$$
{\bf A_0}= (L_{p^0_1,1},...,L_{p_d^0,1}),\;\;\;\;\;{\bf A_1}= (L_{p^1_1,1},...,L_{p_d^1,1}),
$$
and
$$
{\bf B_0}= N_{(q_1^0,...,q_d^0),(\infty,...,\infty)},\;\;\;\;\;{\bf B_1}= N_{(q_1^1,...,q_d^1),(\infty,...,\infty)}.
$$
Taking into account   \cite[Theorem 1]{Nursultanov2000}, we obtain
$$
\left( N_{(q_1^0,...,q_d^0),(\infty,...,\infty)}, N_{(q_1^1,...,q_d^1),(\infty,...,\infty)}\right)_{\bf {\bar \theta},r}\hookrightarrow N_{\bf q,r},
$$
where $
\frac1{\bf q}=\frac{1-{\bar\theta}}{\bf {\bf {\bf
q_0}}}+\frac{\bar\theta}{\bf q_1},$
\, ${\bar \theta}=(\theta_1,...,\theta_d)$ with $0<\theta_i<1.$
Finally,    \cite[Theorem 2]{Nursultanov2000} implies that
$$
L_{\bf p,r} \hookrightarrow
\left( {\bf A_0,A_1}\right)_{\bf {\bar \theta},r}.
$$
Combining the above estimates, we arrive at
(\ref{e13-}).
\end{proof}

The next result is an extension of the known Hardy--Littlewood--Paley inequality for the classical Lorentz spaces,
$$
\big\|\widehat{f}\,\big\|_{L_{p',q}}\lesssim \big\|f\big\|_{L_{p,q}},\quad 1<p<2,\quad0<q \le \infty,
$$
and it is interesting in its own right.
For $\bp=(p_1,...,p_d)$, let $\bp'=(p_1',...,p_d')$.

\begin{thm}\label{t70}
	If $1<\bp<\infty$, $0<\bq \leq \infty$ and $f\in L_{\bp,\bq}$, then  $\mathfrak{F}_N f\in {N_{\bp',\bq}}$  and there holds
	\begin{equation}\label{e13}
		\big\|\mathfrak{F}_N f\,\big\|_{N_{\bp',\bq}}\lesssim \big\|f\big\|_{L_{\bp,\bq}}.
	\end{equation}
	In particular, for $1<p<\infty$,
	\begin{equation}\label{e12}
		  \left(\int_0^\infty \ldots \int_0^\infty (t_1\ldots t_d)^{p-2} \left|\overline{\mathfrak{F}_N f}( t_1,...,t_d;M) \right|^p \, dt\right)^{1/p}
\lesssim \|f\|_{L_p}.
	\end{equation}
Moreover,  the constants $C$ do not depend on $N$.
\end{thm}

\begin{proof} We estimate the $N_{\bp',\infty}$-norm of $\widehat{f}$ as follows:
	\begin{align*}
		\left\|\mathfrak{F}_N{f}\right\|_{N_{\bp',\infty}} &\leq \sup_{I \in M} |I_1|^{-1/p_1}\ldots |I_d|^{-1/p_d} \left| \int_{I_1} \ldots\int_{I_d}\mathfrak{F}_N{f} (\xi_1,...,\xi_d) \, d \xi \right| \\
		&\leq \sup_{I \in M} \int_{Q_N} |f(x)|\prod_{i=1}^d \left(|I_i|^{-1/p_i} \left| \int_{I_i} e^{-\imath \xi_i x_i} \, d \xi_i \right|\right) \, dx\\
&\leq \sup_{I \in M} \int_{\mathbb R^d} |f(x)|\prod_{i=1}^d \left(|I_i|^{-1/p_i} \left| \int_{I_i} e^{-\imath \xi_i x_i} \, d \xi_i \right|\right) \, dx.
	\end{align*}
	Let
	\[
		\varphi_i(x)= |I_i|^{-1/p_i} \left| \int_{I_i} e^{-\imath \xi_i\cdot x_i} \, d \xi_i \right|, \qquad i=1,\ldots,n.
	\]
	If $I_i=[a_i,b_i]$, then
	\begin{align*}
		\left| \int_{I_i} e^{-\imath \xi_i\cdot x_i} \, d \xi_i \right| = \frac{2}{|x_i|} \big| \sin((b_i-a_i)x_i/2)\big| \leq 2 \min \left(\frac{b_i-a_i}{2}, \frac{1}{|x_i|}\right)
	\end{align*}
	and so
	\[
		\varphi_i^{*}(t_i) \leq \frac{4}{(b_i-a_i)^{1/p_i}} \min \left(\dfrac{b_i-a_i}{2}, \frac{1}{t_i}\right) \leq 2t_i^{-1/p_i'}.
	\]
	Using inequality (\ref{hlp})
 $n$-times, we conclude that
	\begin{align*}
		\left\|\mathfrak{F}_N{f}\right\|_{N_{\bp',\infty}} &\leq 2^d \int\limits_{0}^{\infty} t_d^{-1/p'_d} \left(\ldots\int\limits_{0}^{\infty} t_1^{-1/p'_1} f^{*_1}(t_1,\cdot,...,\cdot)\,\frac{dt_1}{t_1}\right)^{*_d}_{t_d} \,\frac{dt_d}{t_d} \\
		&=2^d \left\|\ldots\|f \right\|_{L_{p_1,1}}\ldots \|_{L_{p_d,1}} = 2^d \left\|f \right\|_{(L_{p_1,1},...,L_{p_d,1})} \n.
	\end{align*}
	Using this inequality for $1<\bp_0<\bp<\bp_1<\infty$ and Lemma \ref{t1-}, we derive  that
	\begin{align*}
		\left\|\mathfrak{F}_N{f}\right\|_{N_{\bp',\bq}} &\lesssim \left\|f\right\|_{L_{\bp,\bq}}.
	\end{align*}
	 Setting $p=p_i=q_i,\; i=1,...,n$ in \eqref{e13}, we  immediately obtain  inequality \eqref{e12}.
\end{proof}

\begin{lem}\label{l2}
	For a locally integrable function $f$ and for $I=I_1 \times...\times I_d \in M$ with $|I_i| \leq t_i$, we have
	\[
		\frac{1}{t_1\cdots t_d} \left| \int_{I} f(x)\, dx\right| \leq C_d\, \overline{f}\big( \frac{t_1}2,\cdots, \frac{t_d}2;M).
	\]
\end{lem}

\begin{proof}
Let $I=I_1\times...\times I_d=I_1\times I'$, where   $I'\subset \mathbb R^{n-1}$. Consider a corresponding rectangular parallelepiped $J_1\times I'$, where
 $J_1$ is defined as follows.

If $|I_1|\geq\frac{t_1}{2},$ then define $J_1:=I_1$. In other cases, let $I_1=[a,b]$, $|I_1|=b-a<\frac{t_1}{2}$ and  $Q_1=[a,\,a+t_1]$. If
$$
\left|\int_{I_1\times I'} f(x)dx\right|\leq 2\left|\int_{Q_1\times I'}f(x)dx\right|,
$$
then set   $J_1:=Q_1$. If
$$
\left|\int_{I_1\times I'} f(x)dx\right|>2\left|\int_{Q_1\times I'} f(x)dx\right|,
$$
then set $J_1:=[b,a+t_1]$. Hence,
\begin{eqnarray*}
\left|\int_{J_1\times I'}f(x)dx\right|&=&
\left|\int_{Q_1\times I'}f(x)dx-\int_{I_1\times I'}f(x)dx\right|\\&\geq&\left|\int_{I_1\times I'}f(x)dx\right|-\left|\int_{Q_1\times I'}f(x)dx\right|\geq\frac12\left|\int_{I_1\times I'} f(x)dx\right|.
\end{eqnarray*}
Taking into account that $t \geq |J_1| \geq \frac {t_1}2,$ we derive that
$$
\left|\int_{I_1\times I'}f(x)dx\right|\leq2\left|\int_{J_1\times I'}f(x)dx\right|.
$$
 In this way, in $n$ steps, we get the parallelepiped
 $J=J_1\times ...\times J_d$ such that   $\frac{t_i}2\leq |J_i|\leq t_i$  and
$$
\left|\int_{I}f(x)dx\right|\leq2^d\left|\int_{J}f(x)dx\right|.
$$
Finally, we have

\[
		\frac{1}{t_1\cdots t_d} \left| \int_{I} f(x)\, dx\right| \leq 2^d\frac{1}{t_1\cdots t_d} \left| \int_{J} f(x)\, dx\right|\leq 2^d\, \overline{f}\big( \frac{t_1}2,\cdots, \frac{t_d}2;M).
\]
\end{proof}

\bigskip
\section{Proof of Theorems \ref{t2} and \ref{t3}}\label{s4}

First we obtain the following

\begin{lem}\label{l3}
For a { locally integrable}
 function $f$ satisfying
$$
\sum_{m_d=1}^\infty\ldots\sum_{m_1=1}^\infty \overline{f}(2^{m_1},\ldots,2^{m_d};M)<\infty,
$$
for any $k\in \mathbb Z^d$ and $\varepsilon\in E$, there holds
$$
\sup_{
{}^{2^{k_i}\leq |t_i|\leq2^{k_i+1}}_{\quad i=1,\cdots,n}
}|H_\varepsilon f(t)|
\leq2^d\sum_{m_d=k_d}^\infty\ldots\sum_{m_1=k_1}^\infty \overline{f}(2^{m_1-1},\ldots,2^{m_d-1};M).
$$
\end{lem}
\begin{proof}
For  $t\in \mathbb R^d$, the operator  $H_\varepsilon$ can be represented as
$$
H_\varepsilon f(t)=\int_0^\infty\ldots \int_0^\infty f(t_1\varepsilon_1+x_1\sgn t_1,\ldots,t_d\varepsilon_d+x_d\sgn t_d)\prod_{i=1}^d\psi_{i,\varepsilon_i}(x_i)dx,
$$
with
\[
	\psi_{i,0}(x_i)= \left\{
           \begin{array}{ll}
             \frac{1}{t_i}, & \hbox{if $0 \leq  x_i \leq |t_i|$,} \\
                          0, & \hbox{if $x_i >|t_i|$,}
           \end{array}
         \right.
\quad\mbox{and}\quad\psi_{i,1}(x_i)=\frac1{|t_i|+x_i}.
\]
Let $2^{k_i} \leq |t_i|<2^{k_i+1}, \; i=1,\ldots,n$, then
\begin{align*}
	H_\varepsilon f(t) &=\sum_{m_d=k_d}^\infty \ldots\sum_{m_1=k_1}^\infty\int_{2^{m_d}-2^{k_d}}^{2^{m_d+1}-2^{k_d}}\ldots \int_{2^{m_1}-2^{k_1}}^{2^{m_1+1}-2^{k_1}} \\
	& \qquad f(t_1\varepsilon_1+x_1\sgn t_1,\ldots,t_d\varepsilon_d+x_d\sgn t_d) \prod_{i=1}^d\psi_{i,\varepsilon_i}(x_i)dx.
\end{align*}
Then the mean value theorem gives
$$
H_\varepsilon f(t)=\sum_{m_d=k_d}^\infty\ldots\sum_{m_1=k_1}^\infty\prod_{i=1}^d\psi_{i,\varepsilon_i}(2^{m_i}-2^{k_i})\int_{I_m}f(x)dx,
$$
where  $I_m=I_{m_1}\times\ldots \times I_{m_d}$ is a parallelepiped such that  $|I_{m_i}|\leq 2^{m_i}$. Since
\[
	\psi_{i,\epsilon_i}(2^{m_i}-2^{k_i}) \leq 2^{-m_i} \qquad (i=1,\ldots,n),
\]
 Lemma \ref{l2} completes the proof:
\begin{eqnarray*}
|H_\varepsilon f(t)|&=&\left|\sum_{m_d=k_d}^\infty\ldots\sum_{m_1=k_1}^\infty\prod_{i=1}^d\psi_{i,\varepsilon_i}(2^{m_i}-2^{k_i})\int_{I_m}f(x)dx\right|
\\
&\leq&\sum_{m_d=k_d}^\infty\ldots\sum_{m_1=k_1}^\infty 2^{-m_1-\ldots-m_d}   \left|\int_{I_m}f(x)dx\right|
\\
&\leq&2^d\sum_{m_d=k_d}^\infty\ldots\sum_{m_1=k_1}^\infty \overline{f}(2^{m_1-1},...,2^{m_d-1};M).
\end{eqnarray*}
\end{proof}

\begin{lem}\label{l5}
Let $0<p<\infty$ and $0< q\leq\infty,$.
 Then for  any {  Cauchy sequence }
 $\{f_m\}_m$ from $N_{p,q}$ and any $\varepsilon\in E$ there exists the limit
$$
\lim_{m\to\infty}(H_\varepsilon f_m)(t), \qquad t\in\mathbb R^d.
$$
\end{lem}
\begin{proof}
Let  $\{f_d\}_{n=1}^\infty$ be a Cauchy sequence in $N_{p,q}.$
Let  $\varepsilon>0,$ there exists $N$ such that for $m>N$ and $r\in {\mathbb N}$
one has
\begin{equation}\label{2.2.1}
\|f_{m}-f_{m+r}\|_{N_{p,q}}<\varepsilon.
\end{equation}
Then, by $N_{p,q}\hookrightarrow N_{p,\infty}$, we derive that
$$
\sup_{e\in M}\frac1{|e|^{\frac1{p'}}}\left|\int_e\left(f_m-f_{m+r}\right)d\mu\right| \lesssim \varepsilon, \;\;\;\;\;m>N.
$$
Thus, for any $e\in M$,  $
\left\{\int_e f_m d\mu\right\}_{m=1}^\infty
$ is the
Cauchy sequence, which implies that there exists
\begin{equation}\label{2.2.2}
\lim_{n\rightarrow\infty}\int_{e}f_md\mu.
\end{equation}
Since   $f_m\in N_{p,q}$, Lemma  \ref{l0} yields  that
$$
\sum_{m_d=1}^\infty\ldots\sum_{m_1=1}^\infty \overline{f}(2^{m_1},...,2^{m_d};M)<\infty.
$$

Let  $t\in\mathbb R^d$ and $2^{k_i}\leq t_i<2^{k_i+1}$. By Lemma  \ref{l3} there is $N$ such that
$$
|\left(H_\varepsilon( f_m-f_{m+r})\right)(t)|
\leq 2^{n}\sum_{m_d=k_d-1}^N\ldots\sum_{m_1=k_1-1}^N \overline{(f_m-f_{m+r})}( 2^{m_1},...,2^{m_d};M).
$$
Form the definition of $\overline{f}( t_1,...,t_d;M)$, there are  $Q_m\in M$ satisfying
$$
|\left(H_\varepsilon( f_m-f_{m+r})\right)(t)|
\leq 2^{n+1}\sum_{m_d=k_d-1}^N\ldots\sum_{m_1=k_1-1}^N \frac1{|Q_m|}\left|\int_{Q_m}f_m(x)-f_{m+r}(x)dx\right|.
$$
Finally, taking into account that  the limit \eqref{2.2.2} exists, the sequence  $\left\{\left(H_\varepsilon f_m\right)(t)\right\}_m$
is a Cauchy sequence and, therefore, convergent.
\end{proof}

\begin{proof}[Proof of Theorem \ref{t2}]

Denote $\frac1p=\alpha+\frac1r$. Then we have that $1<p\leq r$ and $\beta=\frac1{p'}-\frac1q$. Due to the embedding  $L_{p,r}\subset L_{p,q}$, taking into account  $\frac1p-\frac1r\geq 0$ and condition \eqref{a1}, we have
\begin{eqnarray}
\|f\|_{L_{p,q}}\lesssim \|f\|_{L_{p,r}}&=&\left(\int_{0}^\infty \ldots\int_{0}^\infty \left((|t_1| \ldots|t_d|)^{\frac1p-\frac1r} f^{*_1...*_d}(t)\right)^{r}\,dt\right)^{1/r}
\nonumber
\\
\label{vvss}
&\lesssim&
 \left(\int_{-\infty}^\infty \ldots\int_{-\infty}^\infty \left((|t_1| \ldots|t_d|)^\alpha |f(t)|\right)^{r}\,dt\right)^{1/r}<\infty.
\end{eqnarray}
Here in the last estimate we have used the Hardy--Littlewood--P\'{o}lya inequality for rear\-rangements
	(see, e.g., \cite[p. 7]{ben})
	\begin{equation}\label{hlp}
 \int_0^{\infty} g^*(t) \frac1{(\varphi^{-1})^*(t)}\, dt
\leq	\int_{-\infty}^{\infty} g(x) \varphi(x)\, dx. 
	\end{equation}
By Theorem \ref{t70} and (\ref{vvss}),
\begin{multline*}
\|\mathfrak{F}_Nf-\mathfrak{F}_{N+r}f\|_{N_{p',q}}
\\
\lesssim\left(\int_{-\infty}^\infty \ldots\int_{-\infty}^\infty \left((|t_1| \ldots|t_d|)^\alpha |(f\chi_{Q_{N+r}\setminus Q_N})(t)|\right)^{r}\,dt\right)^{1/r}\rightarrow 0\quad\mbox{as}\quad N\to +\infty.
\end{multline*}
Thus,   $\left\{\mathfrak{F}_Nf\right\}$ is a Cauchy sequence in $N_{p',q}$, and using Lemma \ref{l5}, there exists
\begin{equation*}
\lim_{N\to+\infty}(H_\varepsilon \mathfrak{F}_N{f})(t).
\end{equation*}
Discretizing the integral 
 yields
\begin{eqnarray*}
		&&\int_{-\infty}^\infty \ldots\int_{-\infty}^\infty \left((|t_1| \ldots|t_d|)^{\beta} \left|\left(H_\varepsilon{\mathfrak{F}_N{f}}\right)(t)\right|\right)^{q}\,dt
\\
&=&	\sum_{\delta\in E}	\int_{0}^\infty \ldots\int_{0}^\infty \left((|t_1| \ldots|t_d|)^{\frac1{p'}-\frac1q} \left|\left(H_\varepsilon{\mathfrak{F}_N{f}}\right)((-1)^{\delta_1}t_1,\ldots,(-1)^{\delta_d}t_d)\right|\right)^{q}\,dt\qquad
\\
&=&	\sum_{\delta\in E}	\sum_{k_d\in \mathbb Z}\ldots\sum_{k_1\in \mathbb Z}\int_{2^{k_d}}^{2^{k_d+1}} \ldots\int_{2^{k_1}}^{2^{k_1+1}} \\
&& \qquad \left((|t_1| \ldots|t_d|)^{\frac1{p'}-\frac1q} \left|\left(H_\varepsilon{\mathfrak{F}_N{f}}\right)((-1)^{\delta_1}t_1,\ldots,(-1)^{\delta_d}t_d)\right|\right)^{q}
 dt.
\end{eqnarray*}
Using Lemmas  \ref{l3}, \ref{l0}, Theorem  \ref{t70} and \eqref{vvss},
  we continue as follows:
\begin{eqnarray*}
&\lesssim&\sum_{k_d\in \mathbb Z}\ldots\sum_{k_1\in \mathbb Z}\left(2^{\frac{k_1+...+k_d}{p'}}\sum_{m_d=k_d-1}^\infty\ldots\sum_{m_1=k_1-1}^\infty \overline{\mathfrak{F}_N{f}}(2^{m_1},\ldots,2^{m_d};M)\right)^q
\\
&\asymp& \| \mathfrak{F}_N{f} \|_{N_{p',q}}^q
\lesssim
\|f\|_{L_{p,q}}^q
 \\
&\,&\qquad\qquad\quad\lesssim \left(\int_{-\infty}^\infty \ldots\int_{-\infty}^\infty \Big((|t_1| \ldots|t_d|)^\alpha |f(t)|\Big)^{r}\,dt\right)^{q/r}.
 \end{eqnarray*}
 Now the theorem follows from Fatou's lemma.

\end{proof}

{

\begin{proof}[Proof of Theorem \ref{t3}]

Since by \eqref{e14}, $H_{\varepsilon}$ in bounded on $L_s$ for $1<s<\infty$, we have
	\begin{align}\label{e49}
		\langle H_{\varepsilon}f,g \rangle = \langle f,H_{1-\varepsilon}g \rangle
	\end{align}
	for $f \in L_s$ and $g \in \cS$. Suppose now that $f$ satisfies condition \eqref{vsp}. Then, by Remark \ref{r1}, $f \in L_s^{loc}$ for some $1<s<\min(r,2)$ and so ${\mathfrak{F}_Nf} \in L_{s'}$. Thus \eqref{e49} can be applied for ${\mathfrak{F}_Nf}$.

	Let us introduce the weighted Lebesgue space $L_{q,\beta}(\R^d)$ by the norm
	$$
\|f\|_{L_{q,\beta}}:=\left(\int_{-\infty}^\infty\ldots \int_{-\infty}^\infty \left((|t_1|\ldots |t_d|)^\beta \left|{f}(t)\right|\right)^qdt\right)^{1/q}.
$$
Since the Schwartz space $\cS$ is dense in $L_{q',-\beta}$, we have
	\[
		\left\|H_\varepsilon {\mathfrak{F}_Nf}\right\|_{L_{q,\beta}} = \sup_{\left\|g\right\|_{L_{q',-\beta}} \leq 1} \langle H_\varepsilon{\mathfrak{F}_Nf},g \rangle = \sup_{\left\|g\right\|_{L_{q',-\beta}} \leq 1} \langle {\mathfrak{F}_Nf},H_{1-\varepsilon}g \rangle = \sup_{\left\|g\right\|_{L_{q',-\beta}} \leq 1} \langle f \chi_{Q_N}, \widehat{H_{1-\varepsilon} g} \rangle,
	\]
	where $g \in \cS$.
 Now we show that $\widehat{H_{1-\varepsilon} g}=H_{\varepsilon} \widehat{g}$. It is enough to prove this for one dimension and for the operator $H$. Indeed, since $Hg \in L_s$ for a given $1<s \leq 2$,
	\begin{align}\label{e51}
	 	\widehat{H g}(y) &= \lim_{N\to \infty} {\mathfrak{F}_N(Hg)(y)} = \lim_{N\to \infty} \int_{Q_N} \frac{1}{x} \int_{0}^{x} g(z) \, dz \, e^{-ixy}\, dx \n \\
	 	&= \lim_{N\to \infty} \int_{Q_N} \int_{0}^{1} g(\xi x) \, d \xi \, e^{-ixy}\, dx\n \\
	 	&= \lim_{N\to \infty} \int_{0}^{1} \int_{Q_N} g(\xi x)  \, e^{-ixy}\, dx \, d \xi \n\\
	 	&= \int_{0}^{1} \int_{\R} g(\xi x)  \, e^{-ixy}\, dx \, d \xi,
	 \end{align}
	 where the limit denotes the $L_{s'}$-limit. The last equation comes from
	 \begin{align}\label{e56}
	 	\left\|\int_{0}^{1} \left(\int_{\R}-\int_{Q_N}\right) g(\xi x)  \, e^{-ix \cdot}\, dx \, d \xi \right\|_{{s'}} 
	 	 \leq \int_{0}^{1} \left\|\left(\int_{\R}-\int_{Q_N}\right) g(\xi x)  \, e^{-ix \cdot}\, dx \right\|_{{s'}}\, d \xi \to 0
	 \end{align}
	 as $N\to \infty$. { Indeed, by Hausdorff-Young inequality,
\begin{align*}
	\left\|\left(\int_{\R}-\int_{Q_N}\right) g(\xi x)  \, e^{-ix \cdot}\, dx \right\|_{{s'}} = \left\|\widehat{G \chi_{\R \setminus Q_N}} \right\|_{{s'}} \leq \left\|G \chi_{\R \setminus Q_N} \right\|_{s} \to 0
\end{align*}
as $N\to \infty$ and $\xi \in (0,1)$, where $G(x):= g(\xi x)$. On the other hand
\begin{align*}
	\left\|G \chi_{\R \setminus Q_N} \right\|_{s} \leq \left\|G \right\|_{s} = \xi^{-1/q}\left\|g\right\|_{s}
\end{align*}
 which is integrable on $(0,1)$ with respect to $\xi$. Now Lebesgue dominated convergence theorem implies \eqref{e56}.
}
Changing the variables in \eqref{e51}, we get that
	 \begin{align*}
	 	\widehat{H g}(y) &= \int_{0}^{1} \int_{\R} g(s)  \, e^{-iys/\xi}\, ds \, \frac{d \xi}{\xi} = \int_{0}^{1} \widehat{g}(y/\xi) \, \frac{d \xi}{\xi} =B \widehat{g}(y).
	 \end{align*}

Using this and H\"{o}lder's inequality, we arrive at 
\begin{align*}
		\left\|H_\varepsilon{\mathfrak{F}_Nf}\right\|_{L_{q,\beta}}
		&= \sup_{\left\|g\right\|_{L_{q',-\beta}} \leq 1} \langle {f}\chi_{Q_N},H_{\varepsilon}\widehat{g} \rangle \\
		& \leq \sup_{\left\|g\right\|_{L_{q',-\beta}}  \leq 1} \int_{-\infty}^\infty\ldots \int_{-\infty}^\infty \left(|t_1|\ldots |t_d|\right)^{\alpha } \left|{f}(t)\chi_{Q_N}(t)\right| \left(|t_1|\ldots |t_d|\right)^{-\alpha } \left|H_{\varepsilon}\widehat{ g}(t)\right| \,dt \n \\
		&\leq \sup_{\left\|g\right\|_{L_{q',-\beta}}  \leq 1} \left(\int_{-\infty}^\infty \ldots\int_{-\infty}^\infty \left( \left(|t_1|\ldots |t_d|\right)^{\alpha} \left|{f}(t)\chi_{Q_N}(t)\right|\right)^{r}\,dt\right)^{1/r} \n\\
		& \qquad\qquad\qquad\left(\int_{-\infty}^\infty \ldots\int_{-\infty}^\infty \left(\left(|t_1|\ldots |t_d|\right)^{-\alpha} \left|H_{\varepsilon}\widehat{g}(t)\right|\right)^{r'}\,dt\right)^{1/r'}. \n
	\end{align*}
		Now, by Theorem \ref{t2}, there holds
	\[
		\left\|H_\varepsilon{\mathfrak{F}_Nf}\right\|_{L_{q,\beta}}  \lesssim \left(\int_{-\infty}^\infty \ldots\int_{-\infty}^\infty\left( \left(|t_1|\ldots |t_d|\right)^{\alpha} \left|{f}(t)\chi_{Q_N}(t)\right|\right)^{r}\,dt\right)^{1/r}.
	\]
{ Finally, the theorem follows from a density argument.}
\end{proof}

}

\bigskip
\section{Hardy spaces and atoms 
} 
\label{s5}

Let us choose a one-dimensional Schwartz function $\phi$ such that $\int_{\R} \phi \, dx \neq 0$. Then we say that a tempered distribution $f$ is in the product Hardy space $H_p=H_p(\R\times \cdots \times \R)$ $(0<p<\infty)$ if
\[
	\left\|f\right\|_{H_p} := \left\|\sup_{t_1>0,\ldots,t_d>0}	\left|(f * \left(\phi_{t_1} \otimes \cdots\otimes\phi_{t_d} \right)\right|\right\|_p <\infty,
\]
where $*$ denotes the convolution, $\phi_s(y):= s^{-1} \phi(y/s)$ $(s>0,y \in \R)$ and
\[
	\left(\phi_{t_1} \otimes \cdots\otimes\phi_{t_d}\right)(x):= \prod_{j=1}^{d} \phi_{t_j}(x_j), \qquad x \in \R^{d}.
\]
It is known that different Schwartz functions yield equivalent norms. Moreover, $H_p$ is equivalent to $L_p$ for $1<p<\infty$. For more about Hardy spaces see \cite{Stein1971,Grafakos2004}.

By a {\it dyadic interval} we mean one of the form $(k2^{-n},(k+1)2^{-n})$. For each dyadic interval $I$ let $I^r$ $(r\in \N)$ be the dyadic interval for which $I\subset I^r$ and $|I^r| = 2^{r} |I|$. If $R:=I_1\times \cdots \times I_d$ is a dyadic rectangle, then set
$R^r := I_1^r \times \cdots \times I_d^r$.

For each dyadic interval $I$ we define $\overline{I} := \left\{x \in  \R: |x| \in \left(|I|^{-1},\infty\right) \right\}$. Obviously, $I\subset J$ implies $\overline{I} \subset \overline{J}$.
For a dyadic rectangle $R=I_1\times \cdots \times I_d$ let
$\overline{R}=\overline{I_1} \times \cdots \times \overline{I_d}$.
If $F\subset \R^d$ is a measurable set, then let
$$
\overline{F} := \bigcup_{R \subset F, R \ {\rm is \ dyadic}} \overline{R}.
$$
It is clear that $F_1\subset F_2$ implies
$\overline{F_1} \subset \overline{F_2}$.

Let us introduce the concept of simple $p$-atoms.
A function $a\in L_2$ is called a {\it simple $p$-atom} if there exist
$I_i\subset \R$ dyadic intervals, $i=1,\ldots,j$ for some $1\leq j \leq d-1$, such that
\begin{enumerate}
\item []
\begin{enumerate}
\item ${\rm supp} \ a \subset I_1\times \ldots I_j \times A$ for some open bounded set $A \subset \R^{d-j}$,
\item
$$ \left\|a\right\|_2 \leq (|I_1| \cdots |I_j| |A|)^{1/2-1/p},$$
\item
$$
\int_{\R} a(x) x_i^k \, dx_i = \int_A a \, d\lambda =0
$$
for all
$i=1,\ldots,j$, $k=0,\ldots ,N=\lfloor 2/p-3/2 \rfloor $ and almost every fixed $x_1,\ldots ,x_{i-1},x_{i+1},\ldots ,x_d$.
\end{enumerate}
\end{enumerate}
If $j=d-1$, we may suppose that $A=I_d$ is also a dyadic interval.
Of course if $a\in L_2$ satisfies these conditions for another subset of $\{1,\ldots,d\}$ than $\{1,\ldots,j\}$, then it is also called simple $p$-atom.

Although not every function in $H_p$ can be decomposed into simple $p$-atoms \cite{Chang1985},
the following result holds.

\begin{lem}\label{l50}
	Let $\eta$ be a measure on the Lebesgue measurable sets of $\R^{d}$ satisfying
	\begin{equation}\label{e140}
		\eta(\overline {F}) \leq C|F| \qquad \mbox{for all open bounded $F\subset \R^d$}.
	\end{equation}
	Let $0<p\leq 1$, $V:L_q\to L_s$ be a bounded linear operator for some $1 \leq q \leq 2$, $1 \leq s \leq  \infty$ and
	\begin{equation*}
		Tf(t)= \left(\prod_{j=1}^{d}t_j\right)^{i} Vf(t), \qquad t\in \R^{d},i=0,1.
	\end{equation*}
	Suppose that there exist $\eta_1,\ldots,\eta_d>0$ such that for every simple $p$-atom $a$ and for every
$r_1\ldots,r_d \in \P$,
\begin{equation}\label{e150}
	\int_{\left(\R \setminus\overline{I_1^{r_1}}\right) \times\cdots\times \left(\R \setminus\overline{I_j^{r_j}}\right)} \int_{\overline{A}} |Ta|^{p} \, d \eta \lesssim 2^{-\eta_1 r_1} \cdots 2^{-\eta_j r_j},
\end{equation}
where $I_1\times \ldots \times I_j \times A$ is the support of $a$.
If $j=d-1$ and $A=I_d$ is a dyadic interval, then we also assume that
\begin{equation}\label{e160}
	\int_{\left(\R \setminus\overline{I_1^{r_1}}\right) \times\cdots\times \left(\R \setminus\overline{I_{d-1}^{r_{d-1}}}\right)} \int_{\left(\overline{I_d}\right)^c} |Ta|^{p} \, d \eta \lesssim 2^{-\eta_1 r_1} \cdots 2^{-\eta_{d-1} r_{d-1}}.
\end{equation}
If $T$ is bounded from $L_2(\R^d)$ to $L_2(\R^d,\eta)$, then
	\begin{align}\label{e59}
		\|Tf\|_{L_p(\R^d,\eta)} \lesssim \|f\|_{H_{p}}, \qquad f\in H_{p} \cap L_q.
	\end{align}
	If $\lim_{k\to\infty} f_k= f$ in $H_p$-norm implies that $\lim_{k\to\infty} Vf_k=Vf$ in the sense of tempered distributions, then \eqref{e59} holds for all $f\in H_p$.
\end{lem}

Note that $H_{p} \cap L_q$ is dense in $H_{p}$.
We omit the proof because it is exactly the same as those of Theorems 3.6.12 and 1.8.1 in \cite{wk4} (see also \cite{wvil}). The only difference is that we have to apply \eqref{e140}. In \cite{wk4}, we supposed that $\overline{F}=F$ and $\eta$ is the Lebesgue measure $\lambda$. For $d=2$ and $\overline{F}=F$, $\eta=\lambda$, the lemma was shown in Fefferman \cite{Fefferman1986a} in a different version. However, that version does not hold for higher dimensions. For $d \geq 3$, the present lemma is due to the last author \cite{wk4}. Applying Lemma \ref{l50}, we can prove Theorems \ref{t7} and \ref{t6}. Since the proofs are more complicated for $d=3$ than for $d=2$ and they are very similar for $d>3$, we present them for $d=3$, only.

\bigskip

\section{Proofs of Theorems \ref{t7}--\ref{t8}} 
\label{s6}

\begin{proof}[Proof of Theorem \ref{t7}]
	Let us introduce the measure
	\begin{equation}\label{e170}
		\eta(A)= \int_A \prod_{j=1}^{d} t_j^{-2} \, dt, \qquad A \subset \R^{d},
	\end{equation}
	and the operator
	\begin{equation*}
		Tf(t)= \left(\prod_{j=1}^{d}t_j\right) \widehat f(t), \qquad t\in \R^{d}.
	\end{equation*}
	We say that $n=(n_1,\ldots,n_d) \in \N^{d}$ and $m=(m_1,\ldots,m_d) \in \N^{d}$ are incomparable if neither $n \leq m$ nor $m \leq n$ hold. Let us denote by $\cF_n$ $(n \in \Z^{d})$ the set of all dyadic rectangles
	$$
	I=(k_12^{-n_1},(k_1+1)2^{-n_1}) \times \cdots \times (k_d2^{-n_d},(k_d+1)2^{-n_d}), \qquad k \in \N^{d}.
	$$
	For $I \in \cF_n$ $(n \in \Z^{d})$, let
	$$
	I_0=(0,2^{-n_1}) \times \cdots \times (0,2^{-n_d}).
	$$
	Since $F$ is bounded, if $I \subset F$, $I \in \cF_n$ is a dyadic rectangle, then $n_1,\ldots,n_d$ are bounded from below. Thus there are only finitely many dyadic rectangles $I^{(j)} \subset F$, $j=1,\ldots,N$ such that $I^{(j)} \in \cF_{n^{(j)}}$ and $n^{(1)},\ldots,n^{(N)}$ are incomparable vectors. It is easy to see that
	\[
		\overline{I^{(j)}} = (2^{n_1^{(j)}},\infty) \times \cdots \times (2^{n_d^{(j)}},\infty), \qquad j=1,\ldots,N
	\]
	and
	\[
		\overline{F} = \bigcup_{j=1}^{N} \overline{I^{(j)}} = \bigcup_{j=1}^{N} \overline{I_0^{(j)}}.
	\]
	For $I^{(j)} \in \cF_{n^{(j)}}$, the union $\cup_{j=1}^{N} I^{(j)}$ has minimal measure if $I^{(j)}\cap I^{(k)} \neq \emptyset$ for all $j \neq k$, more exactly,
	\[
		\left|\bigcup_{j=1}^{N} I_0^{(j)}\right| \leq  \left|\bigcup_{j=1}^{N} I^{(j)}\right|.
	\]
	Indeed, if $I^{(j)}\cap I^{(k)}\neq \emptyset$, then the set $I_0^{(j)}\cup I_0^{(k)}$ arises from the set $I^{(j)}\cup I^{(k)}$ $(j \neq k)$ by a dyadic translation. By the same dyadic translation, we get $I_0^{(j)}\cap I_0^{(k)}$ from $I^{(j)}\cap I^{(k)}$ and the intersections have equal measures. If $I^{(j)}\cap I^{(k)}=\emptyset$, then the set $I_0^{(j)}\cup I_0^{(k)}$ arises from the set $I^{(j)}\cup I^{(k)}$ $(j \neq k)$ by two dyadic translations. The same holds for more than two dyadic rectangles. So the corresponding set to $I_0^{(j)}\cap I_0^{(k)}$ is counted only once in the measure of the union $\cup_{j=1}^{N} I_0^{(j)}$ and at most once in $\cup_{j=1}^{N} I^{(j)}$. By the substitution $1/t_j=x_j$,
	\begin{align*}
		\eta( \overline{F}) &= \eta\left(\bigcup_{j=1}^{N} \overline{I_0^{(j)}}\right) =\int_{\bigcup_{j=1}^{N} \overline{I_0^{(j)}}} \prod_{j=1}^{d} t_j^{-2} \, dt \\
		&= \int_{\bigcup_{j=1}^{N} I_0^{(j)}} 1 \, dx = \left|\bigcup_{j=1}^{N} I_0^{(j)}\right| \leq  \left|\bigcup_{j=1}^{N} I^{(j)}\right| \leq  |F|,
	\end{align*}
	which is exactly \eqref{e140}.

We will prove  \eqref{e150} and  \eqref{e160} for $d=3$, because the proof is similar for larger $d$ or for $d=2$. Choose a simple $p$-atom $a$ with support $R=I_1\times I_2 \times A$ where we may suppose that $I_1=(0,2^{-K_1})$, $I_2=(0,2^{-K_2})$ $(K_1,K_2\in \Z)$. Note that
	\begin{equation*}
	\int_{\left(\R \setminus\overline{I_1^{r_1}}\right) \times \left(\R \setminus\overline{I_2^{r_2}}\right)} \int_{\overline{A}} |Ta|^{p} \, d \eta = \int_{-2^{K_1-r_1}}^{2^{K_1-r_1}} \int_{-2^{K_2-r_2}}^{2^{K_2-r_2}} \int_{\overline{A}} (|t_1| |t_2| |t_3|)^{p-2} \left|\widehat{a}(t)\right|^{p}\,dt.
	\end{equation*}
	By the definition of the atom,
	\begin{align*}
		\left|\widehat a(t)\right|
		&= \left|\int_{I_1} \int_{I_2} \int_{A} a(x_1,x_2,x_3) e^{-\imath t_1x_1} e^{-\imath t_2x_2} e^{-\imath t_3x_3} \, dx_1 \, dx_2 \, dx_3 \right| \n \\
		&= \Bigg| \int_{I_1} \int_{I_2} \int_{A} a(x_1,x_2,x_3) \n\\
		&\qquad \left(e^{-\imath t_1x_1} - \sum_{j=0}^N \frac{(-\imath t_1x_1)^j}{j!} \right) \left(e^{-\imath t_2x_2} - \sum_{j=0}^N \frac{(-\imath t_2x_2)^j}{j!} \right) e^{-\imath t_3x_3} \, dx_1 \, dx_2 \, dx_3 \Bigg| \n \\
		&\leq \int_{I_1} \int_{I_2} \left|e^{-\imath t_1x_1} - \sum_{j=0}^N \frac{(-\imath t_1x_1)^j}{j!} \right| \left|e^{-\imath t_2x_2} - \sum_{j=0}^N \frac{(-\imath t_2x_2)^j}{j!} \right| \n\\
		&\qquad \left|\int_A a(x_1,x_2,x_3) e^{-\imath t_3x_3} \, dx_3 \right| \, dx_1 \, dx_2,
	\end{align*}
	where $N=\lfloor 2/p-3/2 \rfloor $. Using Taylor's farmula,
	\begin{align}\label{e200}
		\left|\widehat a(t)\right|
		&\leq C \int_{I_1} \int_{I_2} |t_1x_1|^{N+1} |t_2x_2|^{N+1} \left|\int_A a(x_1,x_2,x_3) e^{-\imath t_3x_3} \, dx_3 \right| \, dx_1\, dx_2\\
		&\leq C 2^{-K_1(N+1)} 2^{-K_2(N+1)} |t_1|^{N+1} |t_2|^{N+1} \int_{I_1} \int_{I_2} \left|\int_A a(x_1,x_2,x_3) e^{-\imath t_3x_3} \, dx_3 \right| \, dx_1 \, dx_2. \n
	\end{align}
	Then $Np+2p-1 >0$ and
	\begin{align*}
		\int_{\left(\R \setminus\overline{I_1^{r_1}}\right) \times \left(\R \setminus\overline{I_2^{r_2}}\right)} \int_{\overline{A}} |Ta|^{p} \, d \eta
		&\lesssim 2^{-K_1(N+1)p} 2^{-K_2(N+1)p} \\
		&\qquad \int_{-2^{K_1-r_1}}^{2^{K_1-r_1}} \int_{-2^{K_2-r_2}}^{2^{K_2-r_2}} |t_1|^{(N+1)p+p-2} |t_2|^{(N+1)p+p-2} \\
		&\qquad \int_{\overline{A}} |t_3|^{p-2} \left(\int_{I_1} \int_{I_2} \left|\int_A a(x_1,x_2,x_3) e^{-\imath t_3x_3} \, dx_3 \right| \, dx_1 \, dx_2 \right)^p\,dt \\
		&\lesssim 2^{-r_1(Np+2p-1)} 2^{K_1(p-1)} 2^{-r_2(Np+2p-1)} 2^{K_2(p-1)} \\
		&\qquad \int_{\overline{A}} |t_3|^{p-2} \left(\int_{I_1} \int_{I_2} \left|\int_A a(x_1,x_2,x_3) e^{-\imath t_3x_3} \, dx_3 \right| \, dx_1 \, dx_2 \right)^p \,dt_3.
	\end{align*}
	By H\"{o}lder's inequality,
	\begin{align*}
		\int_{\left(\R \setminus\overline{I_1^{r_1}}\right) \times \left(\R \setminus\overline{I_2^{r_2}}\right)} \int_{\overline{A}} |Ta|^{p} \, d \eta
		&\lesssim 2^{-r_1(Np+2p-1)} 2^{K_1(p-1)} 2^{-r_2(Np+2p-1)} 2^{K_2(p-1)} \\
		&\qquad \left(\int_{\overline{A}} |t_3|^{-2} \, dt_3\right)^{(2-p)/2} \\
		& \qquad \left(\int_{\overline{A}} \left(\int_{I_1} \int_{I_2} \left|\int_A a(x_1,x_2,x_3) e^{-\imath t_3x_3} \, dx_3 \right| \, dx_1 \, dx_2 \right)^2 \,dt_3\right)^{p/2} \\
		&\lesssim 2^{-r_1(Np+2p-1)} 2^{K_1(p-1)} 2^{-r_2(Np+2p-1)} 2^{K_2(p-1)} \eta\left({\overline{A}} \right)^{(2-p)/2} \\
		& \qquad \left(\int_{\overline{A}} \left(\int_{I_1} \int_{I_2} \left|\int_A a(x_1,x_2,x_3) e^{-\imath t_3x_3} \, dx_3 \right| \, dx_1 \, dx_2 \right)^2 \,dt_3\right)^{p/2}.
	\end{align*}
	In the next step, we use H\"{o}lder's inequality and Plancherel's theorem and \eqref{e140} to obtain
	\begin{align*}
		\int_{\left(\R \setminus\overline{I_1^{r_1}}\right) \times \left(\R \setminus\overline{I_2^{r_2}}\right)} \int_{\overline{A}} |Ta|^{p} \, d \eta
		&\lesssim 2^{-r_1(Np+2p-1)} 2^{K_1(p-1)} 2^{-r_2(Np+2p-1)} 2^{K_2(p-1)} |A|^{(2-p)/2} |I_1|^{p/2} |I_2|^{p/2} \\
		& \qquad \left(\int_{\overline{A}} \int_{I_1} \int_{I_2} \left|\int_A a(x_1,x_2,x_3) e^{-\imath t_3x_3} \, dx_3 \right|^{2} \, dx_1 \, dx_2 \,dt_3\right)^{p/2} \\
		&\lesssim 2^{-(r_1+r_2)(Np+2p-1)} 2^{(K_1+K_2)(p/2-1)} |A|^{1-p/2} \\
		& \qquad \left(\int_{\overline{A}} \int_{I_1} \int_{I_2} \left|\int_A a(x_1,x_2,x_3) e^{-\imath t_3x_3} \, dx_3 \right|^{2} \, dx_1 \, dx_2 \,dt_3\right)^{p/2} \\
		&\lesssim 2^{-(r_1+r_2)(Np+2p-1)} 2^{(K_1+K_2)(p/2-1)} |A|^{1-p/2} \\
		&\qquad \left( \int_{I_1} \int_{I_2} \int_A |a(x_1,x_2,x_3)|^{2} \, dx \right)^{p/2}.
	\end{align*}
	Taking into account (ii) of the definition of the atom, we conclude
	\begin{align*}
		\int_{\left(\R \setminus\overline{I_1^{r_1}}\right) \times \left(\R \setminus\overline{I_2^{r_2}}\right)} \int_{\overline{A}} |Ta|^{p} \, d \eta
		&\lesssim 2^{-(r_1+r_2)(Np+2p-1)}.
	\end{align*}
	Since $Np+2p-1 >0$, \eqref{e150} holds.

	To prove \eqref{e160}, we obtain
	\begin{equation*}
	\int_{\left(\R \setminus\overline{I_1^{r_1}}\right) \times \left(\R \setminus\overline{I_2^{r_2}}\right)} \int_{\R \setminus\overline{I_3}} |Ta|^{p} \, d \eta = \int_{-2^{K_1-r_1}}^{2^{K_1-r_1}} \int_{-2^{K_2-r_2}}^{2^{K_2-r_2}} \int_{-2^{K_3}}^{2^{K_3}} (|t_1| |t_2| |t_3|)^{p-2} \left|\widehat{a}(t)\right|^{p}\,dt.
	\end{equation*}
	where $I_3=(0,2^{-K_3})$ $(K_3\in \Z)$. Similarly to \eqref{e200},
	\begin{align}\label{e210}
		\left|\widehat a(t)\right|
		&= \left| \int_{I_1}\int_{I_2} \int_{I_3} a(x_1,x_2,x_3) \prod_{k=1}^{3} \left(e^{-\imath t_kx_k} - \sum_{j=0}^N \frac{(-\imath t_kx_k)^j}{j!} \right) \, dx \right| \n\\
		&\leq C \int_{I_1}\int_{I_2} \int_{I_3} |t_1x_1|^{N+1} |t_2x_2|^{N+1} |t_3x_3|^{N+1} | a(x_1,x_2,x_3)| \, dx \n\\
		&\lesssim |t_1|^{N+1} |t_2|^{N+1}|t_3|^{N+1} 2^{-(K_1+K_2+K_3)(N+1)} |I_1|^{1/2}|I_2|^{1/2} |I_3|^{1/2} \n\\
		&\qquad \left(\int_{I_1}\int_{I_2} \int_{I_3} |a(x_1,x_2,x_3)|^{2} \, dx\right)^{1/2} \n\\
		&\lesssim |t_1|^{N+1} |t_2|^{N+1} |t_3|^{N+1} 2^{-(K_1+K_2+K_3)(N+2-1/p)}.
	\end{align}
	Hence,
	\begin{align*}
		\int_{\left(\R \setminus\overline{I_1^{r_1}}\right) \times \left(\R \setminus\overline{I_2^{r_2}}\right)} \int_{\R \setminus\overline{I_3}} |Ta|^{p} \, d \eta
		&\lesssim 2^{-(K_1+K_2+K_3)(Np+2p-1)} \\
		&\qquad\int_{-2^{K_1-r_1}}^{2^{K_1-r_1}} \int_{-2^{K_2-r_2}}^{2^{K_2-r_2}} \int_{-2^{K_3}}^{2^{K_3}} (|t_1| |t_2| |t_3|)^{(N+1)p+p-2} \,dt \\
		&\lesssim 2^{-(r_1+r_2)(Np+2p-1)}.
	\end{align*}
If $\lim_{k\to\infty} f_k= f$ in $H_p$-norm, then the convergence holds also in the sense of tempered distributions and then $\lim_{k\to\infty} \widehat{f_k}= \widehat{f}$ in the sense of tempered distributions. By Lemma \ref{l50}, this completes the proof of Theorem \ref{t7}. 
\end{proof}

\begin{proof}[Proof of Theorem \ref{t6}]
	The proof is similar but slightly more advanced than that of Theorem \ref{t7}.
	We use the measure defined in \eqref{e170} and introduce the operator
	$$
	Tf(t)= \left(\prod_{j=1}^{d}t_j\right) H\widehat f(t), \qquad t\in \R^{d}.
	$$
	Inequality \eqref{e200} implies that
	\begin{align*}
		\left|H \widehat a(t)\right|
		&\leq C 2^{-K_1(N+1)} 2^{-K_2(N+1)} |t_1|^{N+1} |t_2|^{N+1} |t_3|^{-1} \\
		&\qquad \int_{0}^{|t_3|}\int_{I_1} \int_{I_2} \left|\int_A a(x_1,x_2,x_3) e^{-\imath u_3x_3} \, dx_3 \right| \, dx_1 \, dx_2 \, du_3 \\
		&= C 2^{-K_1(N+1)} 2^{-K_2(N+1)} |t_1|^{N+1} |t_2|^{N+1} \int_{I_1} \int_{I_2} H_3|\cF_3a(x_1,x_2,t_3)| \, dx_1 \, dx_2,	
	\end{align*}
	where
	\[
		\cF_3a(x_1,x_2,t_3):= \int_\R a(x_1,x_2,x_3) e^{-\imath t_3x_3} \, dt_3
	\]
	and
	\[
	H_3f(x_1,x_2,t_3):= \frac{1}{t_3} \int_{0}^{t_3} f(x_1,x_2,u_3)\,du_3 \qquad (t_3\neq 0).
	\]
	Then
	\begin{align*}
		&\int_{\left(\R \setminus\overline{I_1^{r_1}}\right) \times \left(\R \setminus\overline{I_2^{r_2}}\right)} \int_{\overline{A}} |Ta|^{p} \, d \eta \\
		&= \int_{-2^{K_1-r_1}}^{2^{K_1-r_1}} \int_{-2^{K_2-r_2}}^{2^{K_2-r_2}} \int_{\overline{A}} (|t_1| |t_2| |t_3|)^{p-2} \left|H\widehat{a}(t)\right|^{p}\,dt \\
		&\lesssim 2^{-K_1(N+1)p} 2^{-K_2(N+1)p} \int_{-2^{K_1-r_1}}^{2^{K_1-r_1}} \int_{-2^{K_2-r_2}}^{2^{K_2-r_2}} |t_1|^{(N+1)p+p-2} |t_2|^{(N+1)p+p-2} \\
		&\qquad \int_{\overline{A}} |t_3|^{p-2} \left(\int_{I_1} \int_{I_2} H_3|\cF_3a(x_1,x_2,t_3)| \, dx_1 \, dx_2 \right)^p\,dt \\
		&\lesssim 2^{-(r_1+r_2)(Np+2p-1)} 2^{(K_1+K_2)(p-1)} \int_{\overline{A}} |t_3|^{p-2} \left(\int_{I_1} \int_{I_2} H_3|\cF_3a(x_1,x_2,t_3)| \, dx_1 \, dx_2 \right)^p \,dt_3.
	\end{align*}
		By H\"{o}lder's inequality,
	\begin{align}\label{e9}
		&\int_{\left(\R \setminus\overline{I_1^{r_1}}\right) \times \left(\R \setminus\overline{I_2^{r_2}}\right)} \int_{\overline{A}} |Ta|^{p} \, d \eta \n\\
		&\lesssim 2^{-(r_1+r_2)(Np+2p-1)} 2^{(K_1+K_2)(p-1)} \n\\
		&\qquad \left(\int_{\overline{A}} |t_3|^{-2} \, dt_3\right)^{(2-p)/2} \left(\int_{\overline{A}} \left(\int_{I_1} \int_{I_2}  H_3|\cF_3a(x_1,x_2,t_3)| \, dx_1 \, dx_2 \right)^2 \,dt_3\right)^{p/2} \n\\
		&\lesssim 2^{-(r_1+r_2)(Np+2p-1)} 2^{(K_1+K_2)(p-1)} \eta\left({\overline{A}} \right)^{(2-p)/2} \n\\
		& \qquad \left(\int_{\overline{A}} \left(\int_{I_1} \int_{I_2}  H_3|\cF_3a(x_1,x_2,t_3)| \, dx_1 \, dx_2 \right)^2 \,dt_3\right)^{p/2}.
	\end{align}
	Taking into account \eqref{e14} and Plancherel's theorem, we conclude that
	\begin{align*}
		&\int_{\overline{A}} \left(\int_{I_1} \int_{I_2}  H_3|\cF_3a(x_1,x_2,t_3)| \, dx_1 \, dx_2 \right)^2 \,dt_3 \\
		&\qquad \leq |I_1| |I_2| \int_{I_1} \int_{I_2} \int_{\R}  \left(H_3|\cF_3a(x_1,x_2,t_3)|\right)^2   \,dt_3\, dx_1 \, dx_2 \\
		&\qquad \leq C|I_1| |I_2| \int_{I_1} \int_{I_2} \int_{\R}  |\cF_3a(x_1,x_2,t_3)|^2   \,dt_3\, dx_1 \, dx_2 \\
		&\qquad \leq C|I_1| |I_2| \int_{I_1} \int_{I_2} \int_{A}  |a(x_1,x_2,x_3)|^2  \, dx_1 \, dx_2 \,dx_3.
	\end{align*}
	Substituting this into \eqref{e9} and using \eqref{e140}, we can see that
	\begin{align*}
		\int_{\left(\R \setminus\overline{I_1^{r_1}}\right) \times \left(\R \setminus\overline{I_2^{r_2}}\right)} \int_{\overline{A}} |Ta|^{p} \, d \eta
		&\lesssim 2^{-(r_1+r_2)(Np+2p-1)} 2^{(K_1+K_2)(p/2-1)} |A|^{1-p/2} \\
		&\qquad \left( \int_{I_1} \int_{I_2} \int_A |a(x_1,x_2,x_3)|^{2} \, dx \right)^{p/2}\\
		&\lesssim 2^{-(r_1+r_2)(Np+2p-1)}.
	\end{align*}
	
	Using \eqref{e210}, we remark that
	\begin{align*}
		\left|H \widehat a(t)\right|&\lesssim |t_1|^{N+1} |t_2|^{N+1} |t_3|^{N+1} 2^{-(K_1+K_2+K_3)(N+2-1/p)}.
	\end{align*}
	Hence, the estimate
	\begin{align*}
		\int_{\left(\R \setminus\overline{I_1^{r_1}}\right) \times \left(\R \setminus\overline{I_2^{r_2}}\right)} \int_{\R \setminus\overline{I_3}} |Ta|^{p} \, d \eta &= \int_{-2^{K_1-r_1}}^{2^{K_1-r_1}} \int_{-2^{K_2-r_2}}^{2^{K_2-r_2}} \int_{-2^{K_3}}^{2^{K_3}} (|t_1| |t_2| |t_3|)^{p-2} \left|H\widehat{a}(t)\right|^{p}\,dt\\
		&\lesssim 2^{-(K_1+K_2+K_3)(Np+2p-1)} \\
		&\qquad\int_{-2^{K_1-r_1}}^{2^{K_1-r_1}} \int_{-2^{K_2-r_2}}^{2^{K_2-r_2}} \int_{-2^{K_3}}^{2^{K_3}} (|t_1| |t_2| |t_3|)^{(N+1)p+p-2} \,dt \\
		&\lesssim 2^{-(r_1+r_2)(Np+2p-1)}
	\end{align*}
	can be proved as in Theorem \ref{t7}. Since $Vf:=H\widehat f$ is bounded from $L_q$ to $L_{q'}$ for all $1 \leq q \leq 2$, Lemma \ref{l50} finishes the proof.
\end{proof}

Let us denote by $BMO$ the dual space of $H_1$ (see \cite{Chang1985}). Note that this space is different from the usual $BMO(\R^{d})$ space, that is the dual of $H_1(\R^{d})$ (see Feffereman and Stein \cite{Fefferman1972}). Similarly to Theorem \ref{t3}, by a duality argument, one can obtain the following corollary. 


\begin{cor}\label{c1}
If $f \in L_1^{loc}$ and
\begin{equation*}
\sup_{t \in \R^{n}} (|t_1|\ldots |t_d| \left|f(t)\right|)<\infty,
	\end{equation*}
	then, for all $N \in \N$,
\begin{equation}\label{e18----}
		\left\|\mathfrak{F}_Nf \right\|_{BMO} \lesssim \sup_{t \in \R^{n}} (|t_1|\ldots |t_d| \left|f(t)\right|),
	\end{equation}
	where $\mathfrak{F}_Nf$ was defined in \eqref{F}.
\end{cor}

Note that inequality \eqref{e22} implies that 
\[
	\left\|H \mathfrak{F}_Nf \right\|_{BMO} \lesssim \sup_{t \in \R^{n}} (|t_1|\ldots |t_d| \left|f(t)\right|).
\]

Note that the similar result to (\ref{e18----})
 for Walsh-Fourier coefficients of one-dimensional functions
 was proved by Ladhawala \cite{Ladhawala1976}. 
Moreover, for the Fourier series $f(x)\sim\sum_{n=0}^\infty a_d e^{ \imath nx}$ with non-negative coefficients $a_n$
the corresponding result follows from a characterization of  BMO  due to
Fefferman (see \cite{sledd}).


\begin{proof}[Proof of Theorem \ref{t8}]
	We use the original version of Lemma \ref{l50}, with $i=0$, $\eta$ the Lebesgue measure and $\overline{F}=F$.
	 It is easy to see that if $a$ is a simple $p$-atom with support $R$ (a dyadic rectangle), then $\supp Ha \subset R$. This means that the integrals in \eqref{e150} and \eqref{e160} are $0$.
	Since $H$ is bounded on $L_q$ for all $1<q<\infty$, Lemma \ref{l50} completes the proof.
\end{proof}

It is known that  the operator $B$ is not bounded from $BMO$ to $BMO$ (see Golubov \cite{Golubov1997}) but the following weaker result holds true.

\begin{cor}\label{c3}
	If $f \in L_\infty\bigcap \bigcup_{1 \leq q<\infty}L_q$, then
	$$
		\left\|Bf\right\|_{BMO} \leq C \left\|f\right\|_{\infty}.
	$$
\end{cor}

\begin{proof}
	Since $B$ in bounded on $L_q$ for $1 \leq q<\infty$, we have
	\begin{align*}
		\langle Bf,g \rangle = \langle f,Hg \rangle,
	\end{align*}
	where $g \in \cS$. We have by Theorem \ref{t8},
	\begin{align*}
		\left\|Bf\right\|_{BMO} &= \sup_{\left\|g\right\|_{H_1} \leq 1} \langle Bf,g \rangle = \sup_{\left\|g\right\|_{H_1} \leq 1} \langle f,Hg \rangle \\
		&\leq \sup_{\left\|g\right\|_{H_1} \leq 1} \|f\|_\infty \|Hg\|_1 \leq C \|f\|_\infty.
	\end{align*}
\end{proof}

However, the operator $H$ is bounded on $BMO$.

\begin{cor}\label{c7}
	If $f \in BMO\bigcap \bigcup_{1<q \leq \infty}L_q$, then
	$$
		\left\|Hf\right\|_{BMO} \leq C \left\|f\right\|_{BMO}.
	$$
\end{cor}

\begin{proof}
	Inequality \eqref{e19} implies that
	\begin{align*}
		\left\|Hf\right\|_{BMO} &= \sup_{\left\|g\right\|_{H_1} \leq 1} \langle Hf,g \rangle = \sup_{\left\|g\right\|_{H_1} \leq 1} \langle f,Bg \rangle \\
		&\leq \sup_{\left\|g\right\|_{H_1} \leq 1} \|f\|_{BMO} \|Bg\|_{H_1} \leq C \|f\|_{BMO},
	\end{align*}
	where $g \in \cS$. In the second equality we used that $H$ in bounded on $L_q$ for $1<q \leq \infty$.
\end{proof}




\bigskip

\section*{Appendix A.
Carleman-type result for Fourier transform}


\begin{exa}\label{ex}
There is a continuous function $F$ from $\cap_{p>1} L_p(\R)$ such that  $\widehat{F}\notin L_p$, $p<2$ and
$$
\int_{\R} {|x|}^{p-2} {|\widehat{F}(x)|}^p d x = \infty,\qquad  p>2.
$$
\end{exa}
Indeed, define
$$
g(x) = \sum\limits_{n=0}^{\infty}
\frac{{\varepsilon}_n}{\sqrt{n+1}\, {\ln }^2 (n+2)}
\chi_{(n-\frac{1}{2}, n+\frac{1}{2})}(x),
$$
where ${\{{\varepsilon}_n \} }_{n=0}^{\infty}$ is the Rudin--Shapiro sequence.

First, it is easy to see that $g \in L_2({\R})$, \, $g
\notin L_p ({\R}),$ \,$ 1<p<2,$ and
$
\int_{\R} {|x|}^{p-2} {|g(x)|}^p dx = \infty,$\; $ p>2.
$
Second,
$$
\int\limits_{-n-\frac{1}{2}}^{n+\frac{1}{2}} e^{-itx} g(x) d x
=\frac{2\sin \frac{t}{2}}{t} \sum\limits_{k=0}^n
\frac{e^{-ikt}{\varepsilon}_k}{\sqrt{k+1}\, {\ln }^2 (k+2)} =:
h(t) f_n (t) .
$$
Applying Abel's transformation,
we obtain
\begin{eqnarray*}
f_n (t) &=& \sum\limits_{k=0}^{n-1} (\frac{1}{\sqrt{k+1}\, {\ln }^2
(k+2)} - \frac{1}{\sqrt{k+2}\, {\ln }^2 (k+3)} )
\sum\limits_{r=0}^k e^{-irt}{\varepsilon}_r \\
&+&\frac{1}{\sqrt{n+1}\, {\ln }^2 (n+2)} \sum\limits_{r=0}^n
e^{-irt}{\varepsilon}_r =:
 \sum\limits_{k=0}^{n-1} a_k P_k (t) + \frac{1}{\sqrt{n+1}\,
{\ln }^2 (n+2)} P_n (t) .
\end{eqnarray*}
Since $|P_k (t)|\leq 5 \sqrt{k+1}$ and  $a_k \leq
\frac{C}{(k+1)\sqrt{k+1}\, {\ln }^2 (k+2)}$, we have that   $f_n\to f$ uniformly, where $f$ is continuous and bounded on ${\R}$, and   $\hat{g}=hf\in L_p ({\R})$ for any $1<p<\infty$.
Finally, we put $F:=\hat{g}$.

\end{document}